\documentclass[bj]{imsart}

\RequirePackage{amsmath, amsthm, amssymb, enumerate, amsfonts, amssymb, mathtools, subcaption, microtype, comment}
\setcitestyle{square, comma, numbers,sort&compress, super}

\RequirePackage{graphicx}
\DeclareMathOperator{\tr}{tr}
\newcommand{\norm}[1]{\left\lVert#1\right\rVert_2}
\newcommand{\snorm}[1]{\left\lVert#1\right\rVert_{\psi_2}}
\newcommand{\fnorm}[1]{\left\lVert#1\right\rVert_{F}}

\newcommand{\ceil}[1]{\lceil {#1} \rceil}

\startlocaldefs
\theoremstyle{plain}

\newtheorem{theorem}{Theorem}[section]
\newtheorem{lemma}[theorem]{Lemma}
\theoremstyle{remark}
\newtheorem{definition}[theorem]{Definition}

\newtheorem{assumption}{Assumption}
\newtheorem*{remark}{Remark}

\newcommand{\uA}       {\mbox{\boldmath$A$}} 
 
\newcommand{\uB}       {\mbox{\boldmath$B$}} 
 
\newcommand{\uC}       {\mbox{\boldmath$C$}}

\newcommand{\uE}       {\mbox{\boldmath$E$}}

\newcommand{\uH}       {\mbox{\boldmath$H$}} 

\newcommand{\uI}       {\mbox{\boldmath$I$}}

\newcommand{\uM}       {\mbox{\boldmath$M$}}

\newcommand{\uO}       {\mbox{\boldmath$O$}}
\newcommand{\uP}       {\mbox{\boldmath$P$}}

\newcommand{\uQ}       {\mbox{\boldmath$Q$}}

\newcommand{\uS}       {\mbox{\boldmath$S$}}

\newcommand{\uU}       {\mbox{\boldmath$U$}} 
 
\newcommand{\uV}       {\mbox{\boldmath$V$}}

\newcommand{\uW}       {\mbox{\boldmath$W$}}  
\newcommand{\uWn}       {\mbox{\boldmath$W_n$}}
  
\newcommand{\uX}       {\mbox{\boldmath$X$}}
\newcommand{\uXn}       {\mbox{\boldmath$X_n$}}
 
\newcommand{\uY}       {\mbox{\boldmath$Y$}}
\newcommand{\uYn}       {\mbox{\boldmath$Y_n$}}

\newcommand{\uZ}       {\mbox{\boldmath$Z$}}  
  
\newcommand{\uiota}             {\mbox{\boldmath$\uiota$}}

\newcommand{\uDelta}            {\mbox{\boldmath$\Delta$}}

\newcommand{\uLambda}           {\mbox{\boldmath$\Lambda$}}

\newcommand{\uSigma}            {\mbox{\boldmath$\Sigma$}}

\newcommand{\uPsi}              {\mbox{\boldmath$\Psi$}}

\newcommand{\uzero}            {\mbox{\boldmath$0$}}
\newcommand{\uone}               {\mbox{\boldmath$1$}}
\newcommand{\BLS}{\hat{\uB}_{LS}}

\endlocaldefs

\begin{document}

\begin{frontmatter}
\title{Posterior consistency in multi-response regression models with non-informative priors for the error covariance matrix in growing dimensions}
\runtitle{Posterior consistency of covariance matrix with non-informative priors}

\begin{aug}
\author[A]{\fnms{Partha}~\snm{ Sarkar$^\ast$}\ead[label=e1]{sarkarpartha@ufl.edu}}
\author[A]{\fnms{Kshitij}~\snm{Khare}\ead[label=e2]{kdkhare@stat.ufl.edu}}
\author[A]{\fnms{Malay}~\snm{Ghosh}\ead[label=e3]{ghoshm@ufl.edu}}
\address[A]{Department of Statistics,
University of Florida\printead[presep={,\ }]{e1,e2,e3}}
\end{aug}

\begin{abstract}
The Inverse-Wishart (IW) distribution is a standard and popular choice of priors for covariance matrices and has attractive properties such as conditional conjugacy. However, the IW family of priors has crucial drawbacks, including the lack of effective choices for non-informative priors. Several classes of priors for covariance matrices that alleviate these drawbacks, while preserving computational tractability, have been proposed in the literature. These priors can be obtained through appropriate scale mixtures of IW priors. However, in the era of increasing dimensionality, the posterior consistency of models that incorporate such priors has not been investigated. We address this issue for the multi-response regression setting ($q$ responses, $n$ samples) under a wide variety of IW scale mixture priors for the error covariance matrix. Posterior consistency and contraction rates for both the regression coefficient matrix and the error covariance matrix are established in the ``large $q$, large $n$'' setting under mild assumptions on the true data-generating covariance matrix and relevant hyperparameters. In particular, the number of responses $q_n$ is allowed to grow with $n$, but with $q_n = o(n)$. Also, some results related to the inconsistency of the posterior distribution and posterior mean for $q_n/n \to \gamma$, where $\gamma \in (0,\infty)$ are provided. 
\end{abstract}


\end{frontmatter}

\section{Introduction}\label{sec1}
Covariance matrix estimation arises in multivariate problems which include multivariate normal models and multi-response regression models. Bayesian estimation of a covariance matrix requires a prior for the covariance matrix. The inverse Wishart (IW) prior has long been used \cite{tom, gelman_book} for the normal covariance matrix with the appealing feature of its conjugacy. However, these priors suffer from two major drawbacks. First, the uncertainty for all variance parameters is controlled by a single degree of freedom parameter which provides no flexibility to incorporate different degrees of prior knowledge to different variance components. Second, commonly used noninformative versions of the IW prior can still impose strong conditions on the diagonal elements of the covariance matrix. For a $q$-dimensional covariance matrix, the default $IW(q+1,\uI)$ prior (a common choice of non-informative IW prior) implies an Inverse-Gamma$(1, 1/2)$ marginal prior for the variances. This marginal prior assigns extremely low density in the region near zero which causes a bias toward larger variances even when the true variance is small. As extensively discussed in references \cite{huang, mulder, zava}, the utilization of an Inverse-Wishart prior on $\uSigma$ leads to substantial interdependence between its diagonal and off-diagonal elements, rendering it unsuitable for scenarios requiring a non-informative prior. Specifically, \cite{tokuda} rigorously illustrates that when $R$ signifies the correlation matrix derived from the covariance matrix $\uSigma\sim IW(q+1,\uI)$, the anticipated effective dependence, denoted by $E[1-|R|^{\frac{1}{q}}]$, approaches 1 as $q$ tends to infinity.

In the one-dimensional case, \cite{gelman2006prior} argued very strongly against the use of vague inverse gamma prior for variances, primarily on the second ground. As an alternative, he recommends the use of uniform and half-t priors for variances. Noting that the half-t prior is a scale mixture of an Inverse-Gamma distribution, \cite{huang} consider an IW distribution with a diagonal scale parameter and assign independent inverse gamma before these diagonal elements. The main advantage of this prior is that one can choose hyperparameters to achieve arbitrarily high non-informativity of all variance and correlation parameters. We refer to this prior as the IG-DSIW (Inverse-Gamma diagonal scale-mixed IW) prior. There are other alternatives as well. For example, keeping the scale parameter of the IW prior as a diagonal matrix \cite{zava} assign independent log-normal and positively truncated normal priors to these diagonal elements. We refer to these priors as the LN-DSIW (Log-normal diagonal scale-mixed IW) prior and TN-DSIW (Truncated-normal diagonal 
scale-mixed IW) prior respectively. On a similar structure, \cite{gelman_book} suggested independent uniform prior distributions for the diagonal scale parameter of the IW prior. Similarly, we will refer to it as U-DSIW (Uniform diagonal scale-mixed IW) prior. For details look at the first subsection of Section \ref{sec3}. (\cite{tokuda}) demonstrated that for a correlation matrix derived from a DSIW prior, the limiting expected effective dependence is reduced to $1-\exp{(-1)}$, which is smaller than that observed with a simple IW prior. From an entirely different perspective, \cite{barnard2000} proposed expressing the covariance matrix in terms of the standard deviations and the correlation coefficients and assigning independent priors for the standard deviations and the correlation coefficients. A key problem with these models is they can often be mostly computationally much slower compared to the scale mixture of IW priors. 

\cite{mulder} considered a matrix-$F$ distribution for the covariance matrix which is a robust version of the inverse Wishart distribution, noting that this distribution reduces to an $F$ distribution in the one-dimensional case. Additionally, matrix-$F$ distributions are also scale-mixed IW distributions with a Wishart distribution for the scale parameter of the IW prior. This representation of the matrix-$F$ distribution leads to additional computational simplicity. There are two other advantages of matrix-$F$ distributions. First, it has ``reciprocity''  properties i.e. when a covariance matrix has a matrix-$F$ distribution then its inverse, the precision matrix also belongs to this family of distributions. Second, matrix-$F$ distribution can be used to construct horseshoe-type priors which are particularly useful for estimating location parameters in a multivariate setting. By choosing hyperparameters carefully one can construct noninformative matrix-$F$ priors. Details are given in Section \ref{sec3}.

While empirical studies in \cite{huang, zava, mulder, gelman_book} have shown promising results, a rigorous study of the asymptotic properties of the posterior resulting from the DSIW and matrix-$F$ priors has not been undertaken. We address this issue in a multi-response regression setting with $q$ responses and $n$ samples. In particular, posterior consistency and contraction rates, both for the regression coefficient matrix and error covariance matrix are established in the ``Large $q$, Large $n$'' setup under some mild and standard assumptions on the data generating system and hyperparameters. This is particularly important considering the increasing dimensionality of many modern datasets. Our consistency result for DSIW prior (Theorem \ref{th_siw}) imposes a very weak restriction on the prior (mixing) densities for the IW diagonal scale entries. These are satisfied by almost all standard densities including the Inverse-gamma (\cite{huang}), Log-normal (\cite{zava}), Truncated-normal (\cite{zava}), and Uniform (\cite{gelman_book}), densities used in previously discussed literature. Thus our results are applicable for a large and flexible class of priors. 

As mentioned before, we establish our posterior consistency results (Theorem \ref{th_siw} and \ref{th_matF}) in a setting that allows $q (= q_n)$ to grow with $n$ but with the restriction that $q_n=o(n)$. We establish the (partial) necessity of this restriction by proving that if $q_n/n \to \gamma $, as $n \to \infty$ where $\gamma \in (0,\infty)$ then the posterior mean will be inconsistent i.e. we will have a serious problem for purposes of estimation (Lemma \ref{incon:siw} and \ref{incon:matF}). In Lemma \ref{incon:whole}, we establish that the entire posterior distribution concentrates on a region that is strictly bounded away from the truth when $q_n/n \to \gamma $, as $n \to \infty$ where $\gamma \in (0,\infty)$. This is not surprising since the priors under consideration do not impose any parameter-reducing restrictions (such as sparsity or low rank) on the covariance matrix. Such restrictions are key to proving consistency results in ultra-high dimensional ``Large $q$, Small $n$'' (or, $q_n>>n$) setting. Adaptation of these priors to the sparsity/ low-rank setting and associated consistency results is beyond the scope of this paper and a basis of future research. Although our main focus is on the parameter $\uSigma$, we also establish a posterior contraction rate for the coefficient matrix $\uB$, with the restriction $p_n = o(n)$ and $q_n = o(n)$ under the spectral norm (Theorem \ref{th_B1}), and with $p_n q_n = o(n)$ under the Frobenius norm (Theorem \ref{th_B2}). This contraction rate is standard and aligns with that of simple least squares estimators for the parameter vector $\uB$ (Lemma \ref{subg4}). It is worth mentioning that the posterior contraction rate for $\uSigma$ aligns with that of the sample covariance matrix of the least squares residuals. This result is embedded in the proofs of Theorems \ref{th_siw} and \ref{th_matF}. 

The rest of the paper is organized as follows. After setting up the basic notation in the next paragraph, we state our basic definitions in Section \ref{sec2}. In Section \ref{sec3}, we discuss the prior and posterior distributions. The main results of this paper are presented in Section \ref{sec4} and Section \ref{sec6}. Section \ref{sec7} addresses inconsistency in estimation and posterior inconsistency in high-dimensional settings. In Section \ref{sec8}, we present the geometric ergodicity proof for the Gibbs sampler under the scale-mixed IW prior. Section \ref{sec9} features numerical simulations. Section \ref{sec5} contains the proofs of the main results from Sections \ref{sec4} and \ref{sec6}. Proofs of certain technical lemmas are provided in the supplemental document (\cite{supp}). The paper concludes with remarks in Section \ref{sec10}.

\textit{Notation.}
Given positive sequences $a_n$ and $b_n$, we shall denote $a_n = O(b_n)$ if there exists a global constant $C$ such that $a_n\leq Cb_n,\;\forall\;n \in \mathbb{N}$ and $a_n = \Omega(b_n)$ if there exist a global constant $C$ such that $a_n\geq Cb_n,\;\forall\;n \in \mathbb{N}$. Also, $a_n = o(b_n)$ to denote the $\lim_{n \to \infty} a_n/b_n = 0$. We will denote the identity matrix of order $n$ by $\uI_n$ and $\uO \in \mathbb{R}^{a \times b}$ denotes the $a \times b$ matrix with all zero entries. Denote $\mathcal{S}^{n-1}$ as a unit Euclidean sphere in $\mathbb{R}^n$. Also, let $\mathbb{P}_q^+$ space of all positive definite matrix of order $q\times q$. For a vector, $x \in \mathbb{R}^n$, $\norm{x}$ denotes its Euclidean norm. Given a metric space $(X, d)$, let $\mathcal{N}(X,\epsilon)$ denote its $\epsilon$-covering number, that is, the minimum number of balls of radius $\epsilon$ needed to cover $X$. We shall use $s_{\min}(\uA)$ and $s_{\max}(\uA)$ to denote the
smallest and largest singular values, respectively for any matrix $\uA$. If $A$ is a symmetric square matrix then $\lambda_{\min}(\uA)$ and $\lambda_{\max}(\uA)$ denote the smallest and largest eigenvalues of $\uA$, respectively. Spectral norm and Frobenius norm are defined in the usual way, with $\norm{\uA}=\sup_{u\in\mathcal{S}^{n-1}}\norm{\uA u}$ and $\fnorm{\uA}=\sqrt{\tr(A^TA)}$.

\section{Preliminaries and model formulation} \label{sec2}
We consider the classical multi-response linear regression setting, viz.
\begin{align}
    \uY=\uX\uB+\uE \label{model1}
\end{align}
where $\uY = (y_1, \cdots , y_n)^{T}$ is an $n \times q$ response matrix obtained from $n$ samples of $q$ continuous response variables, $\uX = (x_1, \cdots, x_n)^T$ is the $n \times p$ design matrix obtained from $n$ samples for $p$ covariates, $\uB \in \mathbb{R}^{p\times q}$
is the coefficient matrix, and $\uE = (\epsilon_1, \cdots , \epsilon_n)^{T}$ is the noise matrix. In particular, $y_i \in \mathbb{R}^q, x_i \in \mathbb{R}^q$ and $\epsilon_i \in \mathbb{R}^q$ for $1 \leq i \leq n$. We assume Gaussianity of the noise, i.e., $\epsilon_1, \cdots, \epsilon_n$ are i.i.d $N_q(0, \uSigma)$. Various {\it working} Bayesian models will be considered in Section \ref{sec3} by combining (\ref{model1}) with priors on the parameters $(\uB, \uSigma)$. 

Next, we recall two key definitions from the literature. 
\begin{definition}{\textbf{(Matrix-normal distribution)}}\label{def1}
A random matrix $\uY$ is said to have a matrix-normal distribution, denoted by $\mathcal{MN}_{n \times q}(\uM,\uU,\uV)$, if $\uY$ has the density function (on the space $\mathbb{R}^{n\times q}$):
\begin{align}
    f(\uY)=\frac{exp\left(-\frac{1}{2}\tr\left[\uV^{-1}(\uY-\uM)^T\uU^{-1}(\uY-\uM)\right]\right)}{(2\pi)^{(np/2)} \vert \uV \vert^{n/2}\vert \uU \vert^{q/2}}\notag
\end{align}
where $\uM\in \mathbb{R}^{n \times q}$, and $\uU$ and $\uV$ are positive definite matrices of dimension $n \times n$  and $q \times q$  corresponding to the variances among the columns and rows of $\uY$, respectively.
\end{definition}
\noindent
Using Definition \ref{def1}, the multi-response linear 
regression model as in (\ref{model1}) can be written simply as 
$\uY\sim \mathcal{MN}_{n \times q}(\uX\uB,\uI_{n},\uSigma)$. 
\begin{definition}{\textbf{(Sub-Gaussian random variable)}}\label{def2} A mean zero random variable $X$ that satisfies $E[\exp(tX)]\leq \exp(t^2k_1^2)$ for all $t \in \mathbb{R}$ and a constant $k_1$ is called a sub-Gaussian random variable. 
\end{definition}

\noindent
If $X$ is a sub-Gaussian random variable then it satisfies $\left( E(|X|^p)\right)^{1/p}\leq k_2\sqrt{p}$ for a constant $k_2$ and one can define its sub-Gaussian norm as $\snorm{X}\coloneqq \sup_{p\geq1}\;p^{-1/2}\left( E(|X|^p)\right)^{1/p}$. A mean zero random vector $X\in \mathbb{R}^q $ is said to be sub-Gaussian if for any $u\in \mathcal{S}^{q-1}$, the random variable $u^{T}X$ is sub-Gaussian. The sub-Gaussian norm of a random vector X is defined as,
\begin{align}
    \snorm{X}\coloneqq \underset{u\in \mathcal{S}^{q-1}}{\sup}\;\snorm{u^{T}X}.\notag
\end{align}

\noindent
See \cite{vershynin} for more details. 

We now specify the {\it true data-generating mechanism}. As previously mentioned, we denote $q$ and $p$ by $q_n$ and $p_n$ to emphasize the fact that $q$ i.e. the number of responses, and $p$ i.e. the number of covariates depends on $n$ respectively and allowed to grow with $n$, so that our results remain relevant in expanding-dimensional settings. We assume that for every $n\geq1$, the observed $q$-dimensional response vectors $y_1^n, y_2^n, \cdots, y_n^n \in \mathbb{R}^{q_n}$ are obtained through 
\begin{align}
      \uYn=\uXn{\uB}_{0n}+\uE_{n},\label{model2}
  \end{align}
where $\uYn :=(y_1^n, . . . , y_n^n)^{T}$, $\uXn$ is the $n \times p_n$ design matrix with the relevant predictor values, and ${\uE}_n\coloneqq(\epsilon_1^n, . . . , \epsilon_n^n)^{T}$ is the noise matrix. We will assume, $\epsilon_1^n, . . . , \epsilon_n^n$ are i.i.d coming from a mean zero sub-Gaussian random variable with $var(\epsilon_1)=\uSigma_{\mathbf{0n}}$. Hence, $\{\uB_{0n}\}_{n\geq1}$ denotes the sequence of true (regression) coefficient matrices and $\{\uSigma_{0n}\}_{n\geq1}$ denotes the sequence of true covariance matrices. For convenience, we denote $\uB_{0n}^{p_n\times q_n}$ as $\uB_{0}$ and $\uSigma_{0n}^{q_n\times q_n}$ as $\uSigma_{0}$ specifically, noting that $\uB_{0}$ depends on $p_n$ and $q_n$ and $\uSigma_{0}$ depends on $q_n$ (and therefore on $n$). Let $\mathbb{P}_{0n}$ denote the probability measure underlying the true model described above and again for convenience we will use $\mathbb{P}_{0}$ instead of $\mathbb{P}_{0n}$. 

Now consider a working Bayesian model based on (\ref{model1}) and a sequence of prior densities $\{\pi_n^*(\uB,\uSigma)\\\}_{n\geq1}$. Let $\{\pi_n^*(\uB,\uSigma|\uYn)\}_{n\geq1}$ denote the sequences of the corresponding joint posterior densities. Analogously, $\{\Pi_n(\cdot)\}_{n\geq1}$ and $\{\Pi_n(\cdot|\uYn)\}_{n\geq1}$ denote the corresponding sequences of joint prior and posterior distributions. We will also use the notation $\pi_n(\cdot)$ and $\Pi_n(\cdot)$ to denote the marginal prior and posterior densities/distributions for $\uB$ and $\uSigma$ when needed. To attain consistent estimation of $\uSigma_{0}$ (denoted as $\uSigma_{0n}$), it is desirable that as the sample size $n$ approaches infinity, the posterior probability assigned outside any $\epsilon$ neighborhood of $\uSigma_{0}$ converges to $0$ in $\mathbb{P}_0$-probability. We, therefore, define a notion of posterior consistency for $\uSigma$ as follows.

\begin{definition}\label{def3}
The sequence of marginal posterior distributions of $\uSigma$ given by $\{\Pi_n(\uSigma|\uYn)\}_{n\geq1}$ is said to be consistent at $\uSigma_{0}$, if for every $\epsilon>0$, $E_{\mathbb{P}_0}[\Pi_n(\norm{\uSigma-\uSigma_{0}} > \epsilon | \uYn)] \to 0$ as $n \to \infty$, where $E_{\mathbb{P}_0}$ denotes an expectation with respect to the joint distribution of $\uYn$ under the true data generating model specified in (\ref{model2}). Additionally, if $E_{\mathbb{P}_0}[\Pi_n(\norm{\uSigma-\uSigma_{0}} > \delta_n | \uYn)] \to 0$ as $n \rightarrow \infty$ for some sequence $\delta_n \rightarrow 0$. Then we denote $\delta_n$ as the contraction rate of $\{\Pi_n(\uSigma|\uYn)\}_{n\geq1}$ around $\uSigma_0$. 
\end{definition}

Although we are mainly interested in the posterior consistency/contraction rate 
for the covariance matrix $\uSigma$, one can similarly define the notion of 
posterior consistency/contraction rate for the regression coefficient matrix $\uB$ (see Section \ref{sec6}). 


\section{Prior and posterior distributions} \label{sec3}

\noindent
Our working Bayesian model combines (\ref{model1}) with a ridge-type matrix normal prior for $\uB$ conditional on $\uSigma$. In particular, we specify 
\begin{equation} \label{beprior}
\uB \mid \uSigma \sim \mathcal{MN}_{p \times q}(0,\;\lambda^{-1}\uI_p,\;\uSigma) 
\end{equation}

\noindent
where $\lambda > 0$ is a user-specified hyperparameter. We will consider two classes of scale-mixed inverse Wishart distributions for the marginal prior of 
$\uSigma$. 

\subsection{The DSIW prior}

\noindent
The diagonal scale-mixed Inverse Wishart (DSIW) prior for $\uSigma$ is based
on the inverse Wishart distribution but adds additional parameters for 
flexibility. The DSIW prior expresses $\uSigma = \uDelta^{1/2} \uQ \uDelta^{1/2}$, where $\uQ$ is positive definite and $\uDelta$ is a diagonal matrix with $i^{th}$ diagonal element equal to $\delta_{i}$. Then, independent priors for $\uQ$ and $\uDelta$ are specified as follows: 
 \begin{align}\label{eq2}
     \uQ \sim IW(\nu+q-1,\;c_{\nu}\uI_q),\;\; \pi(\uDelta) = \prod_{i=1}^q \pi_i (\delta_{i}).  
 \end{align}

 \noindent
 Here $\nu>0$ and $c_{\nu} > 0$ (depends on $\nu$) are user-specified hyperparameters, and $\pi_i(.)$ is a density function with support in a positive real line for every $1 \leq i \leq q$. The DSIW prior can be alternatively represented as 
\begin{align}\label{eq3}
     \uSigma|\uDelta \sim IW(\nu+q-1,\;c_{\nu}\uDelta),\;\;
     \pi(\uDelta) = \prod_{i=1}^q \pi_i (\delta_{i}), 
 \end{align} 
which is easier to handle for theoretical purposes. 

As discussed before, there are several choices of $\pi_i$ (marginal prior on $\delta_i$) prescribed in the literature. \cite{zava} assume independent log-normal priors on $\sqrt{\delta_i}$'s with $c_{\nu}=1$ (LN-DSIW prior). Alternatively, they also suggest using independent truncated normal prior for the $\sqrt{\delta_i}$'s. The induced prior on $\uSigma$ ((TN-DSIW prior) corresponds to the multivariate version of Gelman's folded half-T prior (\cite{gelman2006prior}). \cite{gelman_book} suggest independent uniform prior distributions on $\delta_i$'s with $c_{\nu}=1$ (U-DSIW prior) for a noninformative model. \cite{huang} recommend using independent Gamma prior on $\delta_i$'s with shape parameter $2$ and $c_{\nu}=2\nu$ (IG-DSIW prior), which is a multivariate extension of  Gelman's Half-t priors on standard deviation parameters to achieve arbitrarily high noninformativity. Under this prior correlation parameters have uniform distributions on $(-1, 1)$ when $\nu=2$. This prior is also recommended by \cite{gelman2006data} with $\nu=2$ and $c_{\nu}=1$ to ensure uniform priors on the correlations as in the IW prior but now there is more flexibility on incorporating prior information about the standard deviations. 

We will show in Section \ref{sec4} that consistency for the corresponding DSIW posterior can be obtained with just a mild assumption on the tails of the $\pi_i$s. This assumption is satisfied by all these choices of $\pi_i$ discussed above and provides future researchers with the flexibility of many more choices for $\{\pi_i\}_{i=1}^q$. 

\subsection{Matrix-$F$ prior}
\cite{mulder} propose a matrix-variate generalization of the $F$ distribution, known as the matrix-$F$ distribution, for $\uSigma$. Similar to the univariate $F$ distribution, the matrix-$F$ distribution can be specified via a hierarchical representation as follows: 
\begin{align}\label{eq4}
    \uSigma|\bar{\uDelta} \sim IW(\nu+q-1,\;\bar{\uDelta}),\;\;\bar{\uDelta} {\sim} W(\nu_q^*,\uPsi).
\end{align}

\noindent Here $\nu_q^*>0$ is the user-specified hyperparameter. The key difference in the above specification as compared to the DSIW prior in (\ref{eq3}) is that the scale parameter for the base Inverse Wishart distribution is now a general positive definite matrix. A minimally informative matrix-variate F prior can be obtained by setting $\nu = 1$, $\nu_q^* =q$,  and $\uPsi$ equal to the ``prior guess'' described in \cite{Kass}. The matrix-variate F distribution can also be represented as a Wishart mixture of Wisharts or an inverse Wishart mixture of Wisharts, see \cite{mulder} for further details. 

\subsection{Posterior distributions}
For the working Bayesian models corresponding to both the DSIW prior and the matrix-F priors, the resulting posterior distributions are intractable in the sense that closed-form computations or direct i.i.d. sampling is not feasible. However, the conditional posterior distributions provided in the following two lemmas (proofs provided in the Supplementary Material, \cite{supp}) can be used to construct Gibbs sampling algorithms to generate approximate samples from the resulting posteriors. The conditional/marginal posteriors are also useful in the subsequent theoretical analysis. 
\begin{lemma}\label{Lemma_post1}
For the working Bayesian model which combines (\ref{model1}), (\ref{beprior}) and the DSIW prior in (\ref{eq3}), the posterior density of $(\uB, \uSigma, \uDelta)$ is always proper for any fixed $\nu > 0$. In this case,
\begin{align*}
    &\uB|\uSigma,\uYn \sim \mathcal{MN}_{p \times q}(\tilde{\uB}_n,\;\uX_{\lambda}^{-1},\;\uSigma)\\
    &\uSigma|\uYn,\uDelta\sim IW(\nu+q+n-1,\;\uS_Y+c_{\nu}\uDelta),\\ 
    & \pi_n^{*}(\uDelta \mid \uYn, \uSigma) \propto \prod_{i=1}^q \exp\left\{-c_{\nu}(\Sigma^{-1})_{i}\delta_i/2\right\}  \delta_i^{\frac{\nu+q-1}{2}}\pi_i(\delta_i),
\end{align*}
\noindent where $(\Sigma^{-1})_{i}$ is the $i^{th}$ diagonal element of $\uSigma^{-1}$. Additionally, the marginal posterior density of $\uDelta$ is proportional to 
 \begin{align*}
     \left|\uS_Y+c_{\nu}\uDelta\right|^{-\frac{n+\nu+q-1}{2}}\prod_{i=1}^q \delta_i^{\frac{\nu+q-1}{2}}\pi_i(\delta_i). 
 \end{align*}
 Here, $\uX_{\lambda}=(\uX_n^{T}\uX_n+\lambda \uI_q)$, $\tilde{\uB}_n=\uX_{\lambda}^{-1}\uX_n^{T}\uYn$ and $\uS_Y=\uYn^{T}(\uI_n-\uX_n\uX_{\lambda}^{-1}\uX_n^{T})\uYn.$
\end{lemma}

 \noindent
 The conditional posterior densities of $\uB$ and $\uSigma$ in the above lemma are standard and can be easily sampled from. Note that $\{\delta_i\}_{i=1}^q$ are conditionally independent given $\uYn, \uSigma$. The ease of sampling from the conditional posterior density of $\delta_i$ depends on the choice of the prior density $\pi_i$. For example, if we take $\delta_i \overset{ind}{\sim} Gamma(1/2, A_i^2)$ as suggested by \cite{huang}, then $$\delta_i|\uSigma,\uYn \overset{ind}{\sim} G\left(\frac{\nu+q}{2}, \left(\frac{1}{A_i^2}+\nu(\Sigma^{-1})_{i}\right)^{-1}\right).$$ Even if the conditional posterior of $\delta_i$ is not standard, an efficient accept-reject sampler is available for this univariate density for many choices of $\pi_i$. The next lemma shows that all conditional posterior densities for the matrix-F prior case are standard and easy to sample from. 

\begin{lemma}\label{Lemma_post2}
 For the working Bayesian model which combines (\ref{model1}), (\ref{beprior}) and the matrix-F prior in (\ref{eq4}), the posterior density of $(\uB, \uSigma, \bar{\uDelta})$ is always proper for any fixed $\nu > 0$. In this case,
\begin{align*}
    &\uB|\uSigma,\uYn \sim \mathcal{MN}_{p \times q}(\tilde{\uB}_n,\;\uX_{\lambda}^{-1},\;\uSigma)\\
    &\uSigma|\uYn, \bar{\uDelta}\sim IW(\nu+q+n-1,\;\uS_Y+\bar{\uDelta})\\
    &\bar{\uDelta}| \uSigma,\uYn \sim W(\nu+\nu_q^*+q-1,\;(\uSigma^{-1}+\uPsi^{-1})^{-1}).
\end{align*}
Additionally, the marginal posterior density of $\bar{\uDelta}$ is proportional to 
 \begin{align*}
     \left|\uS_Y+\bar{\uDelta}\right|^{-\frac{n+\nu+q-1}{2}}|\bar{\uDelta}|^{\frac{\nu+\nu_q^*-2}{2}}\exp(-\tr(\uPsi^{-1}\bar{\uDelta})/2).
 \end{align*}
\end{lemma}

\section{Posterior contraction rates for $\uSigma$} \label{sec4}
Before providing our main posterior consistency results, we will state our 
assumptions on the true data-generating model and briefly discuss their implications.
\begin{assumption}\label{as1}
 We assume $q_n=o(n)$.
\end{assumption}

\noindent
We will establish the necessity of this assumption in Section \ref{sec7}. 
As discussed earlier, to obtain consistency results in $q_n = \Omega(n)$ settings, 
priors that impose low-dimensional structure through sparsity, low-rank, etc. are 
needed. 

\begin{assumption}\label{as2} 
There exists $k_{\sigma}\in(0,1]$ such that
  $\uSigma_{0}\in \mathcal{C_{\sigma}}$, where $\mathcal{C_{\sigma}}=\{\uSigma^{q_n\times q_n}| 0<k_{\sigma}\leq\lambda_{min}(\uSigma)\\\leq\lambda_{max}(\uSigma)$ $\leq1/k_{\sigma}<\infty\}$. Here, $k_{\sigma}$ is independent of $n$. We will also assume $\snorm{\uSigma_{0}^{-1/2}\epsilon_1}$ is at most $\sigma_0$, where $\sigma_0 \in \mathbb{R}$ is a fixed number.
\end{assumption}

\noindent
The uniform boundedness of eigenvalues assumption is very standard for consistency results for covariance estimation in both frequentist and Bayesian settings in growing-dimensional setting, see for example, \cite{banerjee1, banerjee2, bickel2008regularized, spectrum, xiang}. \cite{bickel2008regularized} referred to the class $\mathcal{C_{\sigma}}$ as a class of \textit{well-conditioned covariance matrices} and gave many examples of processes that can generate a matrix in this class. The bound on the sub-Gaussian norm (involving $\sigma_0$) is just to ensure that there is no unusual moment behavior for the error distribution. 

\begin{assumption}\label{as3}
  There exists a constant $k$ (not depending on $n$) such that $\pi_i(x)$ decreases in $x$ for $x>k$ for every $1 \leq i \leq q$. 
\end{assumption}
\noindent This assumption is quite mild and covers almost all well-known families of continuous distributions such as truncated normal, half-t distribution, gamma and inverse gamma, beta, Weibull, log-normal, etc. In particular, all the DSIW priors proposed in the current literature satisfy the above assumption for an appropriate choice of $k$. 

\begin{assumption}\label{as4}
 $p_n=O(\sqrt{nq_n})$.
\end{assumption}
\noindent This assumption governs the pace at which the quantity of predictors expands with the variable $n$. For example, this assumption is met when $p_n = O(q_n)$ and when Assumption \ref{as1} is satisfied. It's important to note that we only require this condition to establish a posterior contraction rate for the covariance matrix $\uSigma$. This condition is \textit{not necessary to achieve} a posterior contraction rate for $\uB$, which represents our coefficient matrix. Furthermore, this condition is not a prerequisite for simple posterior consistency for $\uSigma$; rather, its purpose is solely to establish the posterior contraction rate of $\uSigma$. See Remark \ref{remark_post_cont}. 
\begin{assumption}\label{as5}
  There exists $k_x\in(0,1]$ (not depending on $n$) such that $$0<k_x\leq \liminf\limits_{n\rightarrow\infty}\lambda_{min}\left( \frac{\uX_n^{T}\uX_n}{n} \right)\leq\limsup\limits_{n\rightarrow\infty}\lambda_{max} \left( \frac{\uX_n^{T}\uX_n}{n} \right) \leq 1/k_x<\infty.$$ 
\end{assumption}
\noindent Similar assumptions on the design matrix are common in the asymptotic analysis of high dimensional regression models. See for example \cite{bai1, Narisetty:He:2014}.

\begin{assumption}\label{as6}
  $\lambda^{-1}\geq \lambda_{0}\frac{ max(\norm{\uB_0},\norm{\uB_0}^2)}{\sqrt{nq_n}}-\frac{k_x}{n}$, $\lambda_0$ is suitably chosen constant.
\end{assumption}

\noindent
Again, assumptions connecting relevant shrinkage parameters to the norm of the true coefficient vector are typically needed for consistency in high-dimensional regression. As a simple example, consider a simple linear regression model $y = \uX \beta + \epsilon$ with $y \in \mathbb{R}^n$ and $\beta \in \mathbb{R}^p$. The ridge regression estimator with tuning parameter $\lambda$ is 
given by $\hat{\beta}_{ridge} = (\uX^T \uX + \lambda \uI_p)^{-1} \uX^T y$. Suppose the true data generating mechanism is given by $y = \uX \beta_0 + \epsilon$ with 
entries of $\epsilon$ i.i.d. $\mathcal{N}(0, \sigma_0^2)$. Under 
Assumption \ref{as5} on the design matrix, it follows that 
\begin{eqnarray*}
\|\hat{\beta}_{ridge} - \beta_0\|_2
&=& \|(\uX^T \uX + \lambda \uI_p)^{-1} \uX^T (\uX \beta_0 + \epsilon) - \beta_0\|_2\\
&\geq& \|-\lambda (\uX^T \uX + \lambda \uI_p)^{-1} \beta_0\|_2 - \|(\uX^T \uX + \lambda \uI_p)^{-1} \uX^T \epsilon\|_2\\
&\geq& \frac{\lambda}{nk_x + \lambda} \|\beta_0\|_2 - \frac{\epsilon^T {\bf P}_{{\uX}} \epsilon}{nk_x}, 
\end{eqnarray*}

\noindent
Here ${\bf P}_{{\uX}}$ denotes the projection 
matrix into the coulmn space of $\uX$. When $p_n = o(n)$ it can be shown that 
$\epsilon^T {\bf P}_{{\uX}} \epsilon/n$ converges to $0$ in probability. Hence, for consistency of $\hat{\beta}_{ridge}$, it is necessary that 
$\frac{\lambda}{nk_x + \lambda} \|\beta_0\|_2 \rightarrow 0$ as $n \rightarrow \infty$. 

We now state the main results of this paper. The proofs are available in Section \ref{sec5}. 
We first establish the rate of posterior convergence under the DSIW priors on 
$\uSigma$ in Theorem \ref{th_siw}.

\begin{theorem}{\textbf{(Posterior Contraction rate for DSIW prior)}}\label{th_siw}
  Consider a working Bayesian model which combines (\ref{model1}), (\ref{beprior}) with a DSIW prior (\ref{eq3}) on $\uSigma$, and a true data generating mechanism specified in (\ref{model2}). Suppose Assumptions \ref{as1}-\ref{as6} are satisfied. Then for a constant $M>0$ and $\delta_n=\sqrt{q_n/n}$,
  \begin{align*}
      \lim_{n\to\infty} E_{\mathbb{P}_0}\left[\Pi_n(\norm{\uSigma-\uSigma_0} > M\delta_n | \uYn)\right] = 0,
  \end{align*}
  where $E_{\mathbb{P}_0}$ denotes an expectation with respect to the joint distribution of $\uYn$ under the true data generating mechanism. 
\end{theorem}
\begin{remark}
\label{remark_post_cont}
 If the focus is solely on achieving posterior consistency (instead of specifying contraction rate) for the DSIW prior, using a fixed $\epsilon >0$ neighborhood surrounding $\uSigma_0$, then Assumption \ref{as4} can be relaxed to $p_n=o(n)$. This implies that establishing straightforward posterior consistency for the covariance matrix using the DSIW prior requires $p_n=o(n)$ and $q_n=o(n)$ (Assumption \ref{as1}) alongside Assumptions \ref{as2}, \ref{as3}, \ref{as5}, and \ref{as6}.
\end{remark}

\noindent
Even in the basic setting of $n$ i.i.d. $q_n$-dimensional Gaussian observations with common covariance matrix $\uSigma$, the rate of convergence for the sample covariance matrix to its true value can be shown to be $\sqrt{q_n/n}$ under Assumptions \ref{as1} and \ref{as2} (using Lemma \ref{subg2} below, for example, \cite{vershynin, Wainwright}). In the 
 corresponding Bayesian literature, \cite{chaogao} establish a similar contraction rate under Assumptions \ref{as1} and \ref{as2} by employing an $IW(\nu+q_n-1, I_{q_n})$ prior on the covariance matrix. Theorem \ref{th_siw} delineates an equivalent posterior convergence rate within a more intricate and complex multi-response regression framework using scale mixed Inverse-Wishart priors. Theorem \ref{th_siw} shows that adopting a scale mixed Inverse-Wishart prior enables us to harness the aforementioned advantages without any compromise in terms of posterior contraction rate. 

Undoubtedly, there exists a computational cost trade-off, as computations reliant on a scale mixed Inverse-Wishart (IW) prior necessitate an iterative sampling scheme compared to the closed from computation associated with the simple IW prior. This trade-off requires careful consideration by the user, taking into account their resources and specific requirements. Moreover, in Section \ref{sec8}, we establish key theoretical convergence guarantees for the Gibbs sampler associated with a scale mixed IW prior.

We next show in Theorem \ref{th_matF} that the same posterior rate of convergence can be achieved under the matrix-$F$ prior on $\uSigma$ under an additional assumption on a relevant hyperparameter which we state below.

\begin{assumption}\label{as7}
  $\nu_{q}^*=O(q)$.
\end{assumption}

\noindent
The above assumption allows the degrees of freedom of the distribution of the scale parameter $\bar{\uDelta}$ of the base Inverse-Wishart in (\ref{eq4}) to vary with $n$, but stipulates that it cannot be of a larger order than $q_n$. 

\begin{theorem}{\textbf{(Posterior Contraction rate for Matrix-$F$ prior)}}\label{th_matF}
  Consider a Bayesian model which combines (\ref{model1}), (\ref{beprior}) with a Matrix-$F$ prior (\ref{eq4}) on $\uSigma$. Suppose the true data generating mechanism 
  in (\ref{model2}) satisfies Assumptions \ref{as1},\ref{as2} and \ref{as4}-\ref{as7}. Then for a constant $M>0$ and $\delta_n=\sqrt{q_n/n}$,
  \begin{align*}
      \lim_{n\to\infty} E_{\mathbb{P}_0}\left[\Pi_n(\norm{\uSigma-\uSigma_0} > M\delta_n | \uYn)\right] = 0,
  \end{align*}
  where $E_{\mathbb{P}_0}$ denotes an expectation with respect to the joint distribution of $\uYn$ under the true data generating mechanism. 
\end{theorem}

\section{Posterior contraction rates for $\uB$}\label{sec6}

\noindent
In this section, we examine the posterior contraction rate of the regression coefficient matrix $\uB$ under a ridge-type matrix normal prior for $\uB$ given 
$\uSigma$ in (\ref{beprior}), and both the DSIW and matrix-F prior choices for the marginal prior of $\uSigma$. It turns out that the only role the choice of the marginal prior for $\uSigma$ plays in the analysis below is the need for the corresponding (marginal) posterior of $\uSigma$ to be consistent. Since such property has been established for the DSIW and matrix-F priors in Section \ref{sec4}, we present a unified proof for both of these prior choices below. 
Recall the posterior contraction rate for $\uSigma$ in Theorems \ref{th_siw} and \ref{th_matF} is $\delta_n = \sqrt{q_n/n}$. When employing the spectral norm as a distance measure, a contraction rate of $\delta_n^{*}=\max\{\sqrt{p_n/n},\;\sqrt{q_n/n}\}$ can be established for $\uB$. Alternatively, when utilizing the Frobenius norm, the contraction rate becomes $\delta_n^{**}=\sqrt{p_n q_n/n}$. When $q_n$ remains constant (independent of $n$), it becomes evident that the contraction rate under both norms will be $\sqrt{p_n/n}$. This is because the operator norm and Frobenius norm coincide for a vector. Notably, this rate aligns with various optimal rates found in the literature, as demonstrated in works such as those by \cite{song,castillo, wand,cai, zhang}. Furthermore, these rates of posterior contraction align with those of the simple least square estimator for the parameter vector $\uB$ within the framework of the model specified in (\ref{model1}) (refer to Lemma \ref{subg4}). Thus, we observe analogous rates for the entire posterior distribution of $\uB$ as compared to those obtained through a frequentist estimator. Nevertheless, the distinct advantage of our Bayesian approach lies in its natural ability to offer uncertainty quantification.
\begin{theorem}{\textbf{(Posterior Contraction rate for $\uB$ under spectral norm)}}\label{th_B1}
  Consider a Bayesian model which combines (\ref{model1}) with the specified matrix normal 
  prior for $\uB$ (see (\ref{beprior})) given $\uSigma$, and either a DSIW prior or a Matrix-$F$ prior on $\uSigma$. Suppose Assumptions \ref{as1}-\ref{as3},\ref{as5}-\ref{as7} are satisfied together with the assumption that $p_n= o(n)$. Then for a constant $M>0$ and if $\delta_n^{*}=\max\{\sqrt{p_n/n},\;\sqrt{q_n/n}\}$
  \begin{align*}
      \lim_{n\to\infty} E_{\mathbb{P}_0}[\Pi_n(\norm{\uB-\uB_0} > M\delta_n^{*} | \uYn)] = 0,
  \end{align*}
  where $E_{\mathbb{P}_0}$ denotes an expectation with respect to the joint distribution of $\uYn$ under the true data generating mechanism. 
\end{theorem}

\noindent To determine the posterior contraction rate of $\uB$ under the Frobenius norm, we must adjust Assumption \ref{as6} by incorporating $\fnorm{\uB}$, denoted as Assumption \ref{as8}:
\begin{assumption}\label{as8}
  For a suitably chosen constant $\lambda_0$, we have $\lambda^{-1} \geq \lambda_{0}\frac{ \max(\fnorm{\uB_0},\fnorm{\uB_0}^2)}{\sqrt{np_nq_n}}-\frac{k_x}{n}$.
\end{assumption}
\noindent The implication of Assumption \ref{as8} remains consistent with that of Assumption \ref{as6}; the only modification involves substituting the spectral norm with the Frobenius norm. This adjustment aligns with our focus on achieving posterior consistency under the Frobenius norm.
\begin{theorem}{\textbf{(Posterior Contraction rate for $\uB$ under Frobenius norm)}}\label{th_B2}
  Consider a Bayesian model that combines (\ref{model1}) with the specified matrix normal 
  prior for $\uB$ ( see (\ref{beprior})) given $\uSigma$, and either a DSIW prior or a Matrix-$F$ prior on $\uSigma$. Suppose Assumptions \ref{as2},\ref{as3},\ref{as5},\ref{as7},\ref{as8} are satisfied together with the assumption that $p_n q_n= o(n)$. Then for a constant $M>0$ and if $\delta_n^{**}=\sqrt{p_nq_n/n}$
  \begin{align*}
      \lim_{n\to\infty} E_{\mathbb{P}_0}[\Pi_n(\fnorm{\uB-\uB_0} > M\delta_n^{**} | \uYn)] = 0,
  \end{align*}
  where $E_{\mathbb{P}_0}$ denotes an expectation with respect to the joint distribution of $\uYn$ under the true data generating mechanism. 
\end{theorem}
\noindent Note that to determine a posterior contraction rate for the parameter vector $\uB$, we don't need Assumption \ref{as4}. The growth condition required for Theorem \ref{th_B1} is that the number of covariates $p_n$ should be $o(n)$. In contrast, for Theorem \ref{th_B2}, a slightly stronger assumption is needed: $p_n q_n = o(n)$. This distinction arises because the Frobenius norm is larger than the spectral norm. This is because establishing a posterior contraction rate for $\uB$ only requires posterior consistency for the covariance matrix $\uSigma$. The rate at which the posterior distribution contracts for $\uSigma$ does not impact the posterior contraction rate for $\uB$. The proofs for Theorems \ref{th_B1} and \ref{th_B2} can be found in the Section \ref{sec5}.

\section{Inconsistency of posterior distribution and posterior mean when $q_n \neq o(n)$}\label{sec7}

\noindent
A natural question to ask is whether Assumption \ref{as1} i.e. $q_n=o(n)$ can be relaxed to include ultra high-dimensional settings where $q_n$ is of the same order or much larger than $n$. We will show in this section that if $q_n/n\to\gamma\in(0,\infty)$ as $n \to \infty$, then the posterior mean of $\uSigma$ under both the DSIW and matrix-$F$ priors will be an inconsistent estimator of the true covariance matrix even in a much simpler i.i.d. setting. Consider in particular observations $y_1, y_2, \cdots, y_n \in \mathbb{R}^q$ which are assumed to be i.i.d. from a multivariate normal distribution with mean $0$ and covariance matrix $\uSigma$. Suppose, Assumption \ref{as2} is valid for the true underlying covariance matrix, $\uSigma_0$. As a consequence, in an asymptotic regime where $q=q_n$ grows at the same rate as $n$, i.e., $q_n/n\to\gamma\in(0,\infty)$, the spectral norm of the difference of sample covariance matrix and true covariance matrix can be shown to be bounded away from zero, see \cite{bai}. We leverage this result to establish that the posterior mean of $\uSigma$ under both the DSIW and matrix-F priors is inconsistent (Lemmas \ref{incon:siw} and \ref{incon:matF} below). Inconsistency of the posterior mean does not necessarily imply inconsistency of the posterior distribution. From a Bayesian standpoint, one may be interested in the behavior of the entire posterior distribution when $q_n = \Omega(n)$, as opposed to focusing solely on a point estimate like the posterior mean. In the same i.i.d normal setting, it can be shown that when $q_n/n\to\gamma\in(0,\infty)$, the resulting posterior is also not consistent (Lemma \ref{incon:whole}). These results strongly indicate that posterior consistency in $q_n = \Omega(n)$ settings is in general infeasible unless an additional low-dimensional structure is induced in the covariance matrix. First, we will state a technical lemma needed to prove our inconsistency results. 

\begin{lemma}\label{chisq}
 Suppose $q_n$ is a sequence of integers such that $\lim_{n\to\infty}\frac{q_n}{n}=\gamma\in[0,1]$ with $q_n\leq n$ and if $\{S_i\}_{i=1}^n$ is a sequence of random variables such that $nS_i\overset{i.i.d}{\sim}\chi^2_{n}$. Then if we define $Q_{q_n}
 ^{*}=\underset{1\leq i \leq q_n}{\max}S_i$ and $Q_{q_n}
 ^{**}=\underset{1\leq i \leq q_n}{\min}S_i$ , then $Q_{q_n}^{*}$ and $Q_{q_n}^{**} \to 1$ almost surely as $n \to \infty$.
\end{lemma}

\noindent Next, in Lemma \ref{incon:siw} and \ref{incon:matF}, we demonstrate the inconsistency of the posterior mean in the setting where $q_n = \Omega(n)$. The proofs are provided in the Supplementary Material \cite{supp}.
\begin{lemma}\label{incon:siw}
 Consider a Bayesian model where $y_1, y_2, \ldots, y_n$ are i.i.d. Gaussian with mean $0$ and covariance matrix $\uSigma$, and a DSIW prior is placed on $\uSigma$ as described in (\ref{eq3}) with $\nu=1$ and $c_{\nu}=1$. Under Assumption \ref{as3} on $\pi_i(.)'s$ and Assumption \ref{as2} on true underlying covariance matrix $\uSigma_0$, if $q_n/n\to\gamma\in(0,\infty)$, then $\exists\;\omega>0$ such that 
 \begin{align*}
     \liminf_{n\to\infty}{\mathbb{P}_0\left(\norm{\hat{\uSigma}_n-\uSigma_0}>\omega\right)}>0,
 \end{align*}
 where $\hat{\uSigma}_n$ is the posterior mean. However, if $q_n=o(n)$ under the same setup, then
 \begin{align*}
     \norm{\hat{\uSigma}_n-\uSigma_0}\overset{P}{\rightarrow}0.
 \end{align*}

\end{lemma}

\begin{lemma}\label{incon:matF}
 Consider a Bayesian model where $y_1, y_2, \ldots, y_n$ are i.i.d. Gaussian with mean $0$ and covariance matrix $\uSigma$, and a matrix-$F$ prior is placed on $\uSigma$ as described in (\ref{eq4}) with $\nu=1$ and $\nu_{q}^*=q_n$.  Under Assumption \ref{as2} on true underlying covariance matrix $\uSigma_0$, if $q_n/n\to\gamma\in(0,\infty)$, then $\exists\;\omega>0$ such that 
 \begin{align*}
     \liminf_{n\to\infty}{\mathbb{P}_0\left(\norm{\hat{\uSigma}_n-\uSigma_0}>\omega\right)}>0,
 \end{align*}
 where $\hat{\uSigma}_n$ is the posterior mean. However, if $q_n=o(n)$ under the same setup, then
 \begin{align*}
     \norm{\hat{\uSigma}_n-\uSigma_0}\overset{P}{\rightarrow}0.
 \end{align*}
\end{lemma}

\noindent Note that the inconsistency of the posterior mean does not necessarily imply inconsistency of the entire posterior distribution. Thus, in Lemma \ref{incon:whole}, we establish the inconsistency of the posterior distribution under the condition where $q_n = \Omega(n)$. The proof of this lemma is furnished in the Supplementary Material \cite{supp}.

\begin{lemma}\label{incon:whole}
  Consider a Bayesian model where $y_1, y_2, \ldots, y_n$ are i.i.d. Gaussian with a mean of $0$ and a covariance matrix $\uSigma$. A DSIW prior (as described in (\ref{eq3})) or a matrix-$F$ prior (as described in (\ref{eq4})) or a simple $IW(\nu+q-1,\uPsi)$ ( where $\uPsi$ is a non-random  matrix) prior is placed on $\uSigma$. Assuming $\uSigma_0=I_{q_{n}}$, if $q_n/n\to\gamma\in(0,\infty)$, then there exists $\omega>0$ such that 
\begin{align*}
   \liminf_{n\to\infty}  E_{\mathbb{P}_0}[\Pi_n(\norm{\uSigma-\uSigma_{0}} > \omega \mid \uYn)] > 0,
\end{align*}
where $E_{\mathbb{P}_0}$ denotes an expectation with respect to the joint distribution of $\uYn$ under the true data generating model.
\end{lemma}

\noindent In Lemma \ref{incon:whole}, we chose the true covariance matrix value, i.e., $\uSigma_0=I_{q_{n}}$, to simplify our computations. This conclusion of the lemma remains valid for a general $\uSigma_0$ as long as Assumption \ref{as2} holds. While \cite{chaogao} establish posterior contraction rates for a simple IW prior when $q_n=o(n),$ there's a gap in existing literature concerning the inconsistency of the posterior distribution of the covariance matrix, even with a simple IW prior, when $q_n=\Omega(n)$. Lemma \ref{incon:whole} fills this gap and encompasses nearly all conceivable variants of the IW prior or the scaled mixed IW prior found in the literature. The strength of Lemma \ref{incon:whole} lies in the simplicity of its requirement—only the positive definiteness of the scale matrix of the original IW prior is needed. We don't require any additional assumptions on the hyperparameters.

\section{Geometric ergodicity of Gibbs samplers for scale-mixed IW prior}\label{sec8}

\noindent
In this section, we establish the geometric ergodicity of the Gibbs samplers corresponding to a large class of DSIW priors and the matrix-F prior. Geometric ergodicity guarantees that the distribution of the associated Markov chain converges at a geometric rate to the desired posterior distribution. Furthermore, it also implies the existence of an associated Markov chain CLT, which can be leveraged to compute asymptotically valid standard errors for Markov chain-based estimates of posterior quantities, see \cite{Jones:2004, VFJ:2019}.

We will first consider the DSIW prior. Using the conditional posterior densities in Lemma \ref{Lemma_post1} one can easily implement a Gibbs sampling Markov chain $(\uB^{(i)},\uSigma^{(i)},\uDelta^{(i)})_{i=0}^{\infty}$. Suppose, $f_{DSIW} (\cdot \mid \cdot)$ is the Markov transition density (with respect to the Lebesgue measure on $\mathbb{R}_+^q$) corresponding to the $\uDelta$-marginal chain $(\uDelta^{(i)})_{i=0}^{\infty}$. With 
$\uDelta^{(i)}= diag(\delta_1^{(i)},\cdots,\delta_q^{(i)})$, we establish a geometric drift condition using the drift function $V_1(\uDelta)=\sum_{i=1}^{q}(\delta_i+\frac{1}{\delta_i})$ via the following lemma. 
\begin{lemma}\label{ge1}
Suppose the DSIW prior in (\ref{eq3}) satisfies the assumption 
$$
c_2\;x^{\alpha_1-1}\exp\{-\frac{x}{\alpha_2}\} \leq \pi_i(x) \leq c_1\;x^{\alpha_1-1}\\\exp\{-\frac{x}{2\alpha_2}\} \; \forall x > 0
$$
\noindent for some positive constants $\alpha_1$ and $\alpha_2$. Then there exist constant $\lambda_1 \in [0,1)$ and $b_1\geq 0$ such that 
\begin{align*}
    E_{\uDelta \sim f_{DSIW} (\cdot \mid \uDelta^{(0)})}[V_1(\uDelta)\mid V_1(\uDelta^{(0)})]\leq \lambda_1 V_1(\uDelta^{(0)}) +b_1,
\end{align*}
for every $\uDelta^{(0)}\in \mathbb{}{R}_+^q$.
\end{lemma}

\noindent
Since $V_1(\uDelta) < d$ implies $\delta_j \in (1/d, d)$ for every $d > 1$, it follows that the function $V_1$ is {\it unbounded off compact sets} on $\mathbb{R}_+^q$. Also, the form of the conditional densities in Lemma \ref{Lemma_post1} implies that $f_{DSIW} (\uDelta \mid \tilde{\uDelta})$ is continuous in $\tilde{\uDelta}$. A standard argument based on Fatou’s Lemma \cite[Page 127]{meyn1993markov} implies that the $\uDelta$-marginal chain is a Feller chain. Now, a sequential application of Theorem 6.0.1 and Lemma 15.2.8 in \cite{meyn1993markov}, combined with the geometric drift condition established above, implies geometric ergodicity of the $\uDelta$-marginal chain. Note that the DSIW Gibbs sampler based on Lemma \ref{Lemma_post1} has a {\it two-block} structure: it alternates between sampling the block $(\uB, \uSigma)$ given (the current value of) $\uDelta$, and then samples 
$\uDelta$ given (the current value of) $(\uB, \uSigma)$. Hence, the $\uDelta$-marginal chain can be viewed as a de-initializing chain for the Gibbs sampling Markov chain $(\uB^{(i)},\uSigma^{(i)},\uDelta^{(i)})_{i=0}^{\infty}$. Geometric ergodicity of the DSIW Gibbs sampling chain now follows immediately by results in \cite{Roberts:Rosenthal:2001}. Finally, the condition required for the DSIW prior in Lemma \ref{ge1} is relatively mild and encompasses the Inverse Gamma DSIW (IG-DSIW) (\cite{huang}, independent Inverse Gamma priors on $\delta_i$s) prior, and Truncated Normal DSIW (TN-DSIW) prior (\cite{zava}, independent truncated normal priors on $\sqrt{\delta_i}$s) by making appropriate choices for the parameters $\alpha_1$ and $\alpha_2$.

For matrix-$F$ prior, using the conditional posterior densities in Lemma \ref{Lemma_post2} one can easily implement a Gibbs sampling Markov chain $(\uB^{(i)},\uSigma^{(i)},\bar{\uDelta}^{(i)})_{i=0}^{\infty}$. Suppose, $f_{MF} (\cdot \mid \cdot)$ is the Markov transition density (with respect to the Lebesgue measure on $\mathbb{P}_q^+$) corresponding to the $\bar{\uDelta}$-marginal chain $(\bar{\uDelta}^{(i)})_{i=0}^{\infty}$. We establish a geometric drift condition using the drift function $V_2(\bar{\uDelta})=\tr(\bar{\uDelta})+\tr(\bar{\uDelta}^{-1})$ via the following lemma. 
\begin{lemma}\label{ge2}
Consider the matrix-$F$ prior. Then there exist constant $\lambda_2 \in [0,1)$ and $b_2\geq 0$ such that 
\begin{align*}
    E_{\bar{\uDelta} \sim f_{MF} (\cdot \mid \bar{\uDelta}^{(0)})}[V_2(\bar{\uDelta})\mid V_2(\bar{\uDelta}^{(0)})]\leq \lambda_2 V_2(\bar{\uDelta}^{(0)}) +b_2,
\end{align*}
for every $\bar{\uDelta}^{0}\in \mathbb{P}_q^+$.
\end{lemma}

\noindent
Similar arguments as in the DSIW setting above can be used in conjunction with Lemma \ref{ge2} to establish geometric ergodicity of the matrix-$F$ Gibbs sampling chain $(\uB^{(i)},\uSigma^{(i)},\bar{\uDelta}^{(i)})_{i=0}^{\infty}$. The proofs for Lemma \ref{ge1} and \ref{ge2} are available in the supplementary material (\cite{supp}).

\begin{figure}
\begin{subfigure}{0.4\textwidth}
  \centering
  \includegraphics[width=\linewidth]{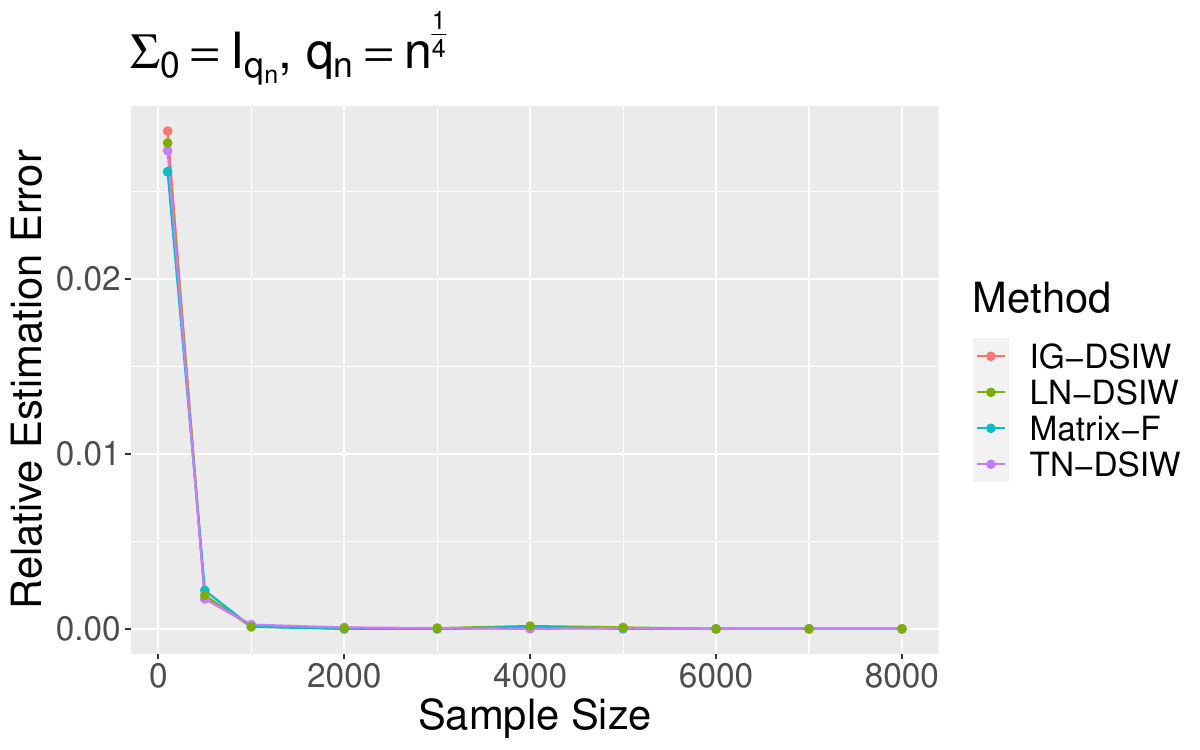}

\end{subfigure}%
\begin{subfigure}{0.4\textwidth}
  \centering
  \includegraphics[width=\linewidth]{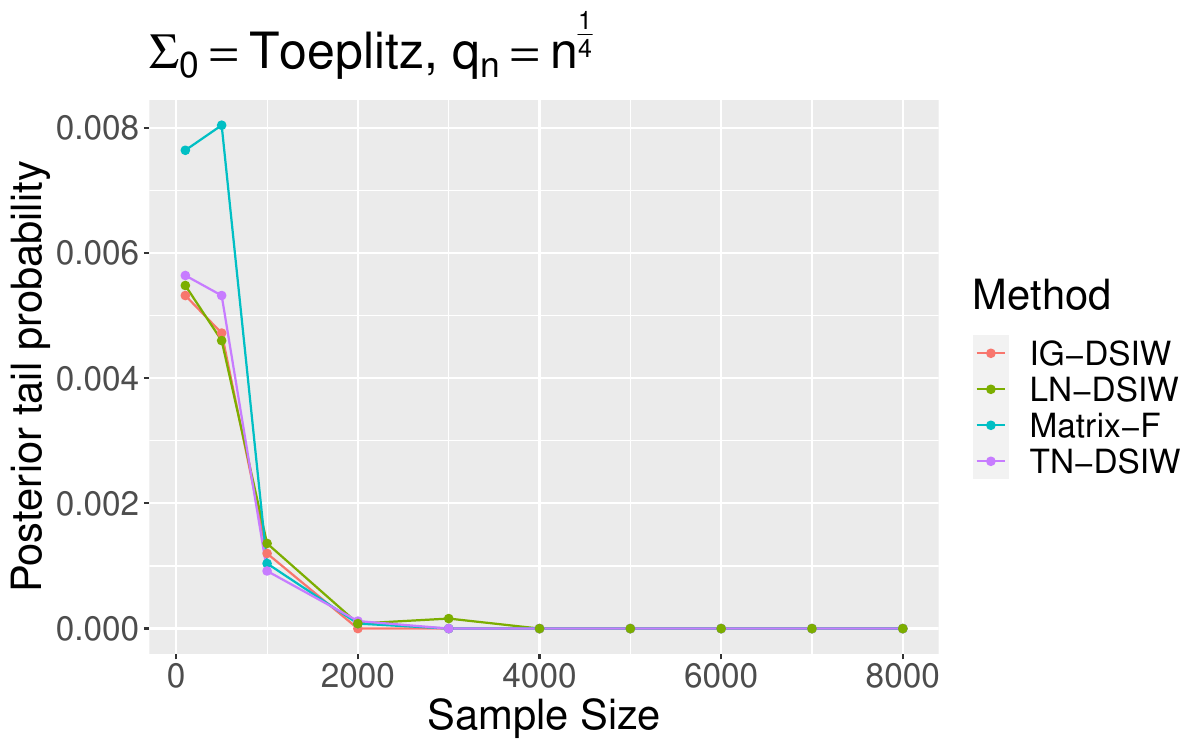}

\end{subfigure}
\begin{subfigure}{0.4\textwidth}
  \centering
  \includegraphics[width=\linewidth]{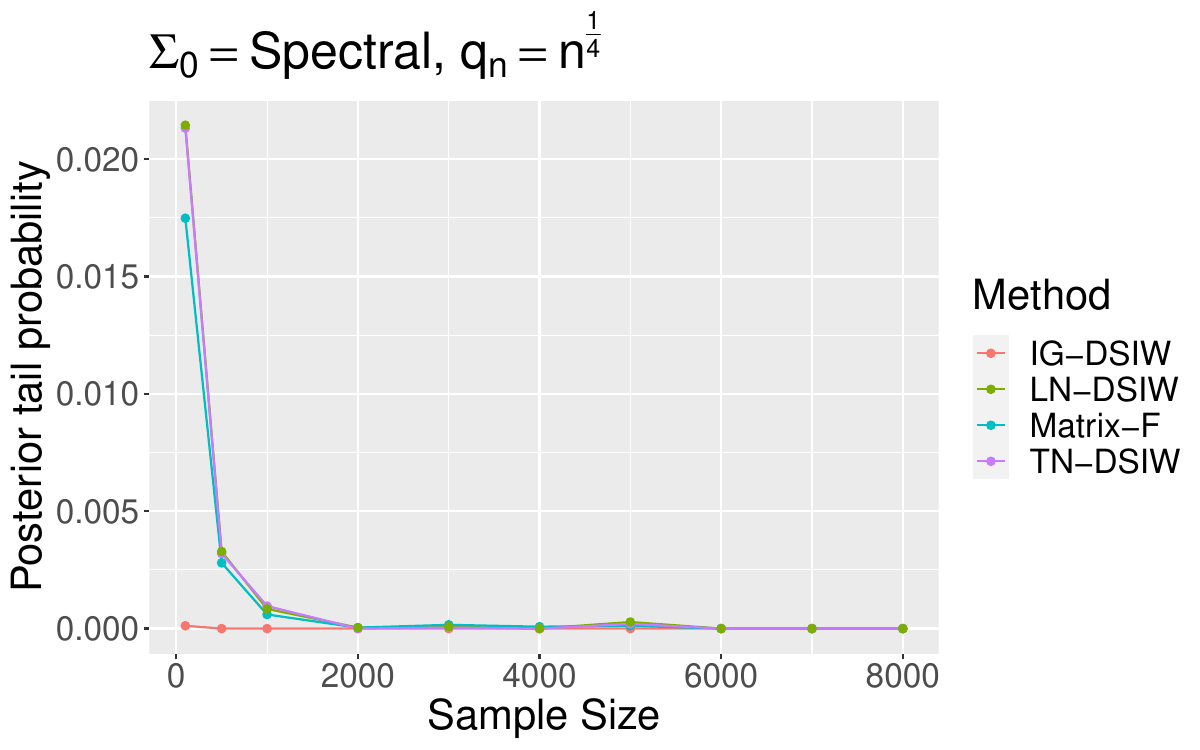}

\end{subfigure}%
\begin{subfigure}{.4\textwidth}
  \centering
  \includegraphics[width=\linewidth]{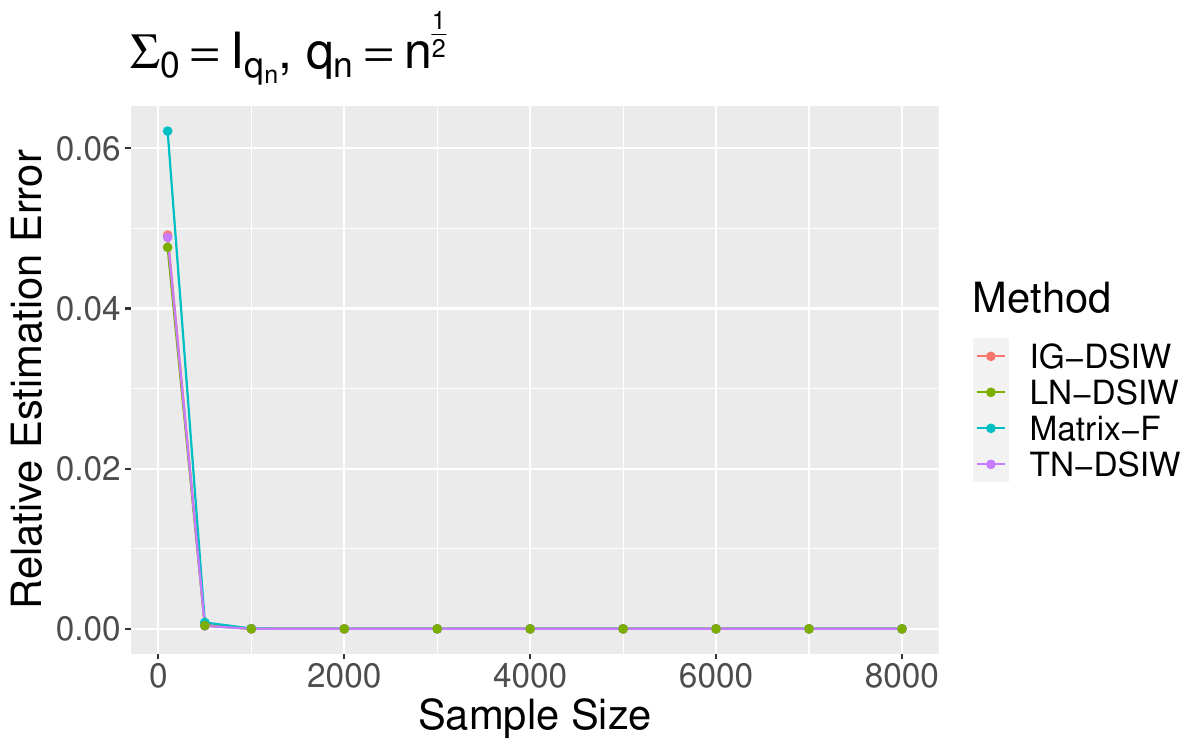}
  
\end{subfigure}
\begin{subfigure}{.4\textwidth}
  \centering
  \includegraphics[width=\linewidth]{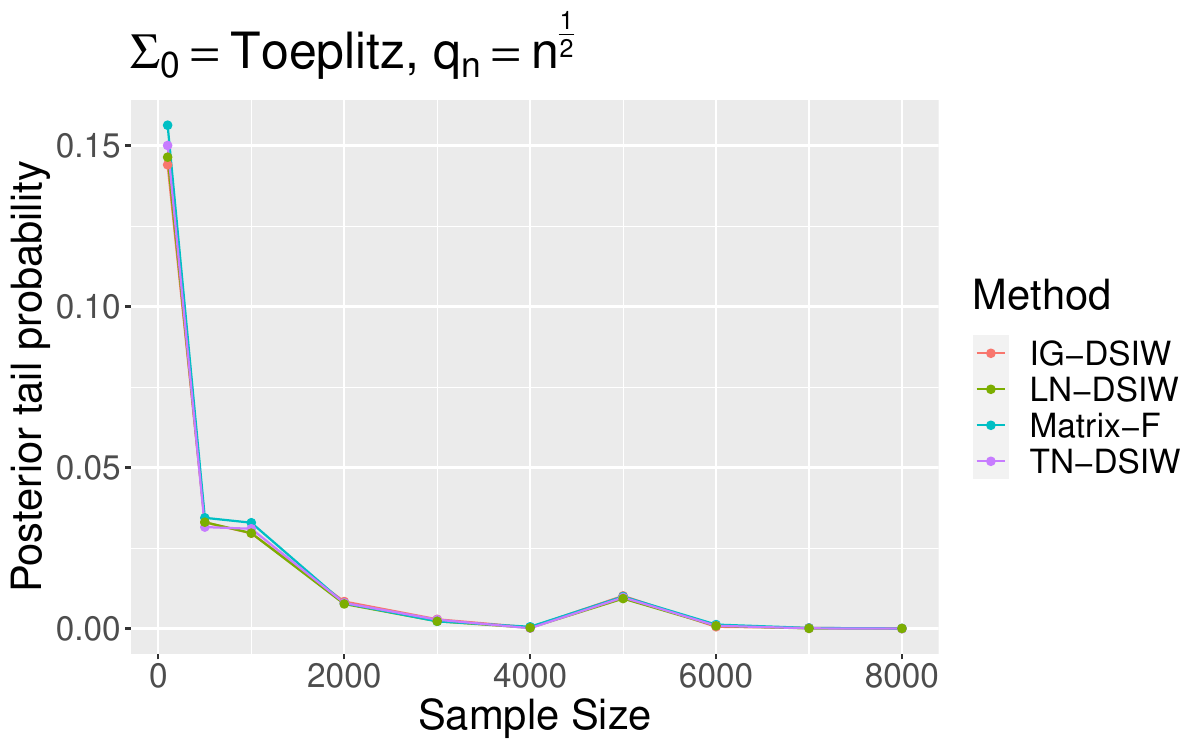}

\end{subfigure}%
\begin{subfigure}{.4\textwidth}
  \centering
  \includegraphics[width=\linewidth]{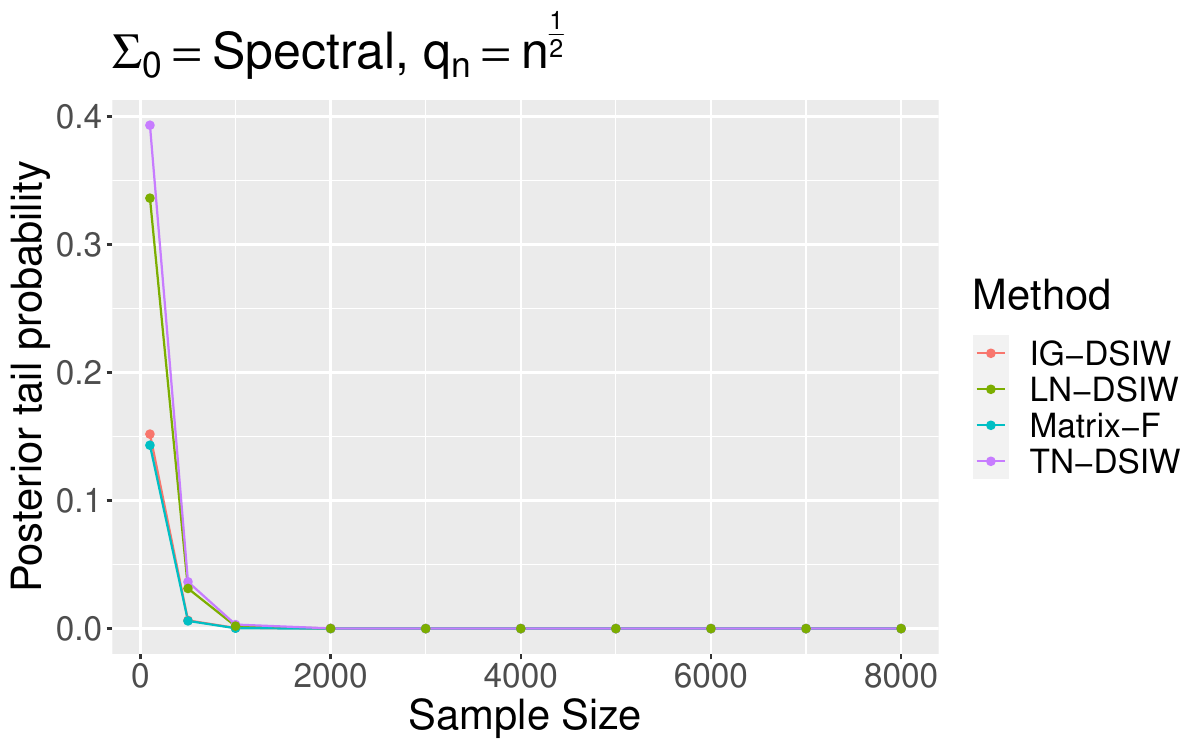}

\end{subfigure}
\begin{subfigure}{.4\textwidth}
  \centering
  \includegraphics[width=\linewidth]{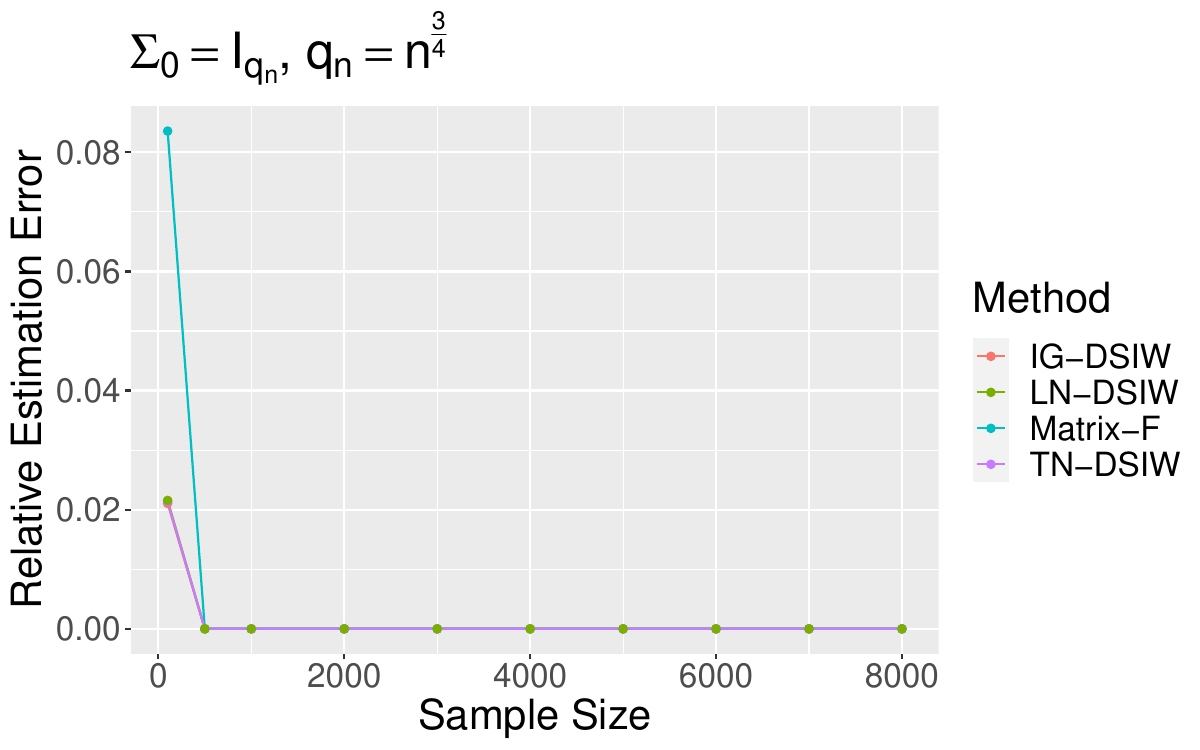}

\end{subfigure}%
\begin{subfigure}{.4\textwidth}
  \centering
  \includegraphics[width=\linewidth]{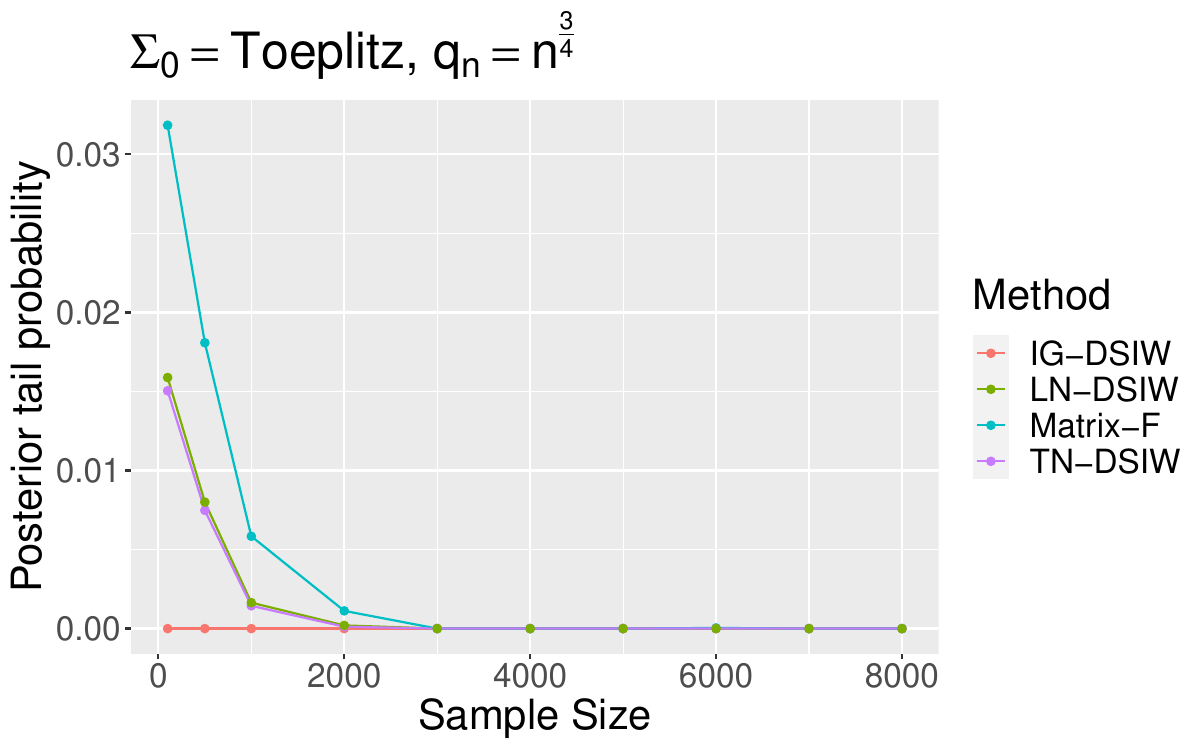}

\end{subfigure}
\begin{subfigure}{0.4\textwidth}
  \centering
  \includegraphics[width=\linewidth]{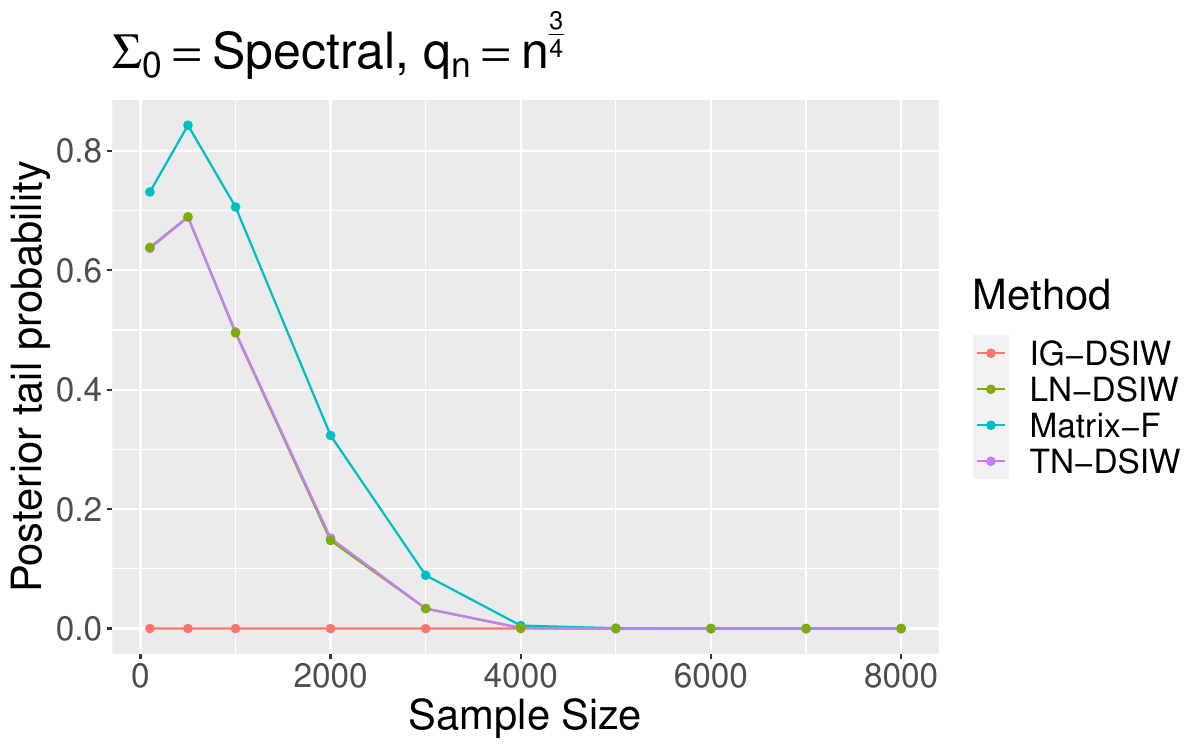}
  
\end{subfigure}
\caption{Plots of posterior tail probability $\Pi_n(\norm{\uSigma-\uSigma_0} > M\sqrt{q_n/n}\mid\uYn)$ vs. sample size $n$ for nine different combinations of $\uSigma_0$ and $q_n$ (with $q_n/n \rightarrow 0$).}
\label{fig:fig1}
\end{figure}

\begin{figure}
\begin{subfigure}{0.4\textwidth}
  \centering
  \includegraphics[width=\linewidth]{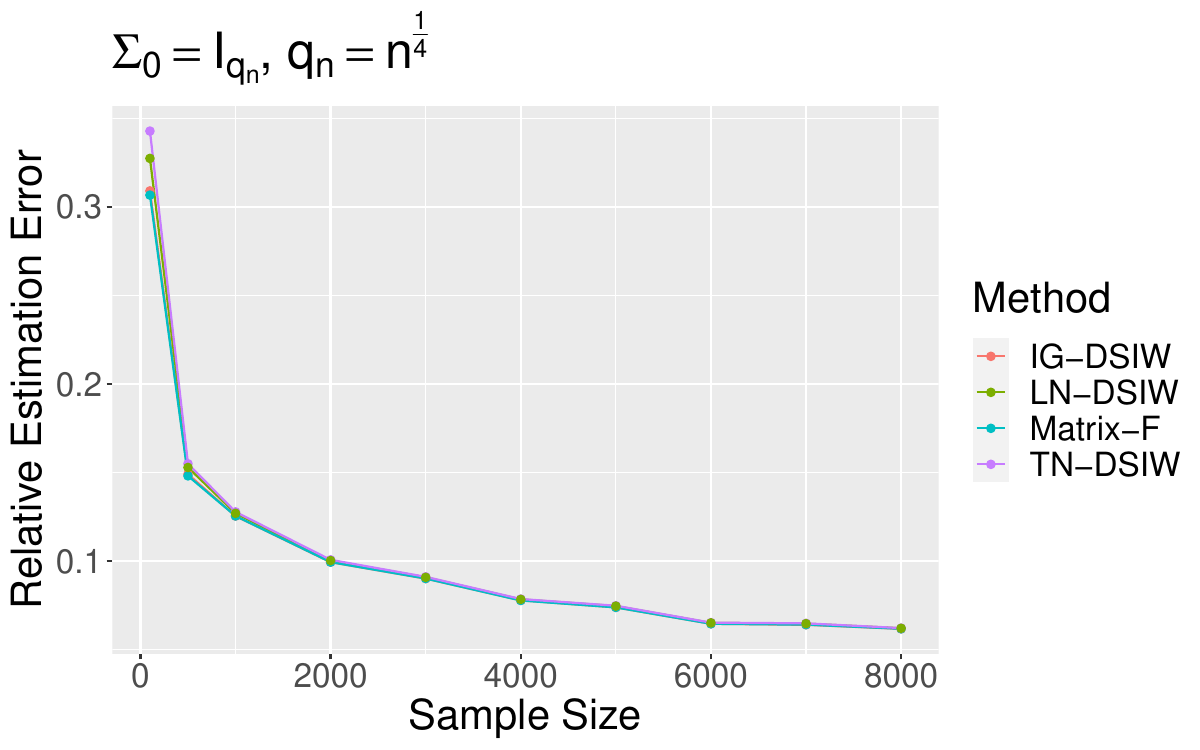}

\end{subfigure}%
\begin{subfigure}{0.4\textwidth}
  \centering
  \includegraphics[width=\linewidth]{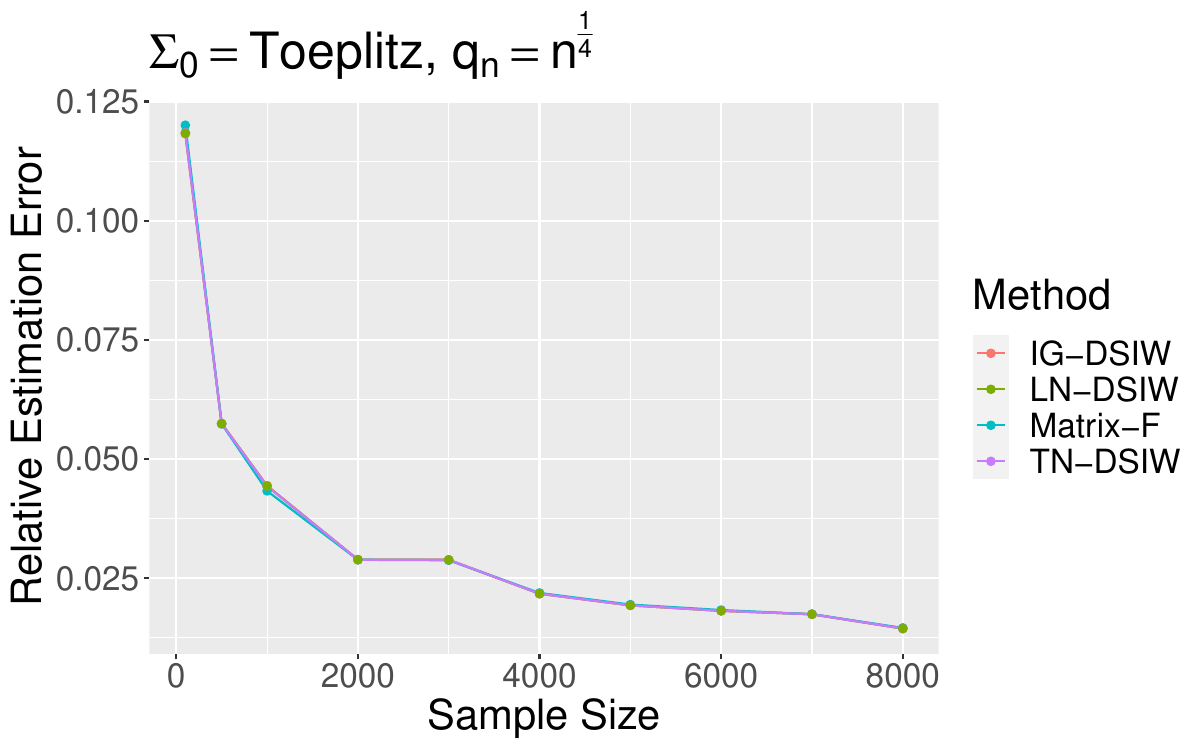}

\end{subfigure}
\begin{subfigure}{0.4\textwidth}
  \centering
  \includegraphics[width=\linewidth]{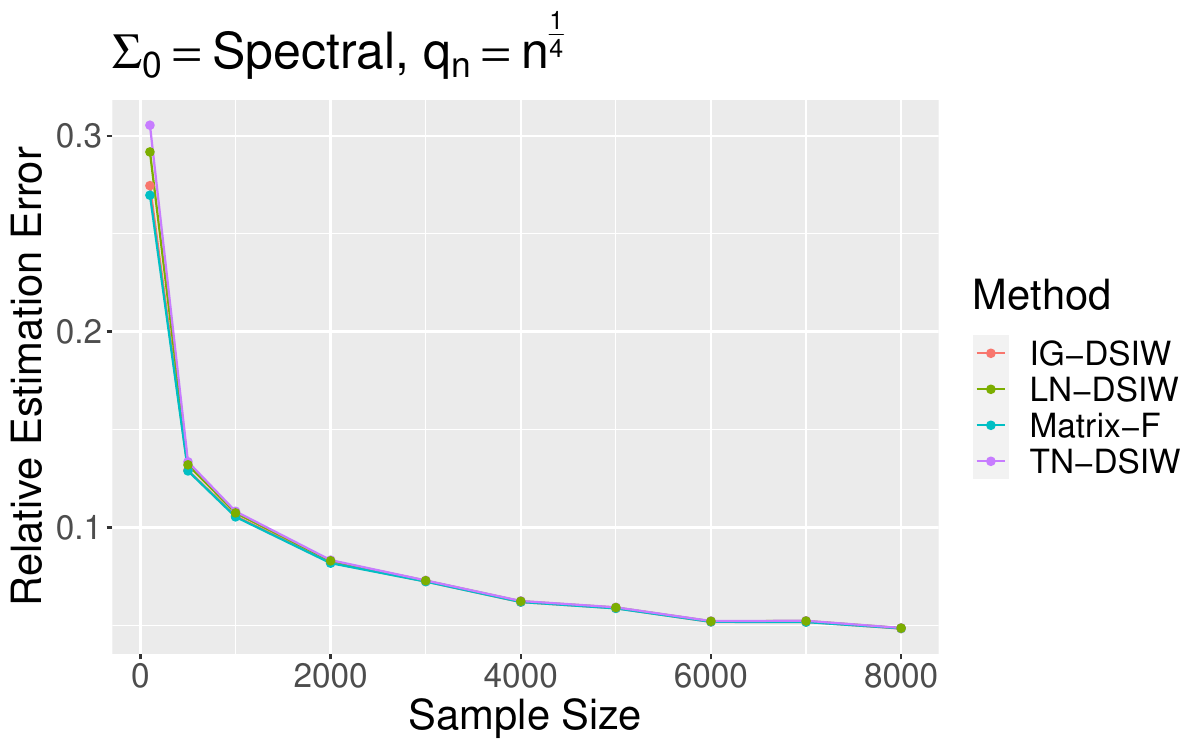}

\end{subfigure}%
\begin{subfigure}{.4\textwidth}
  \centering
  \includegraphics[width=\linewidth]{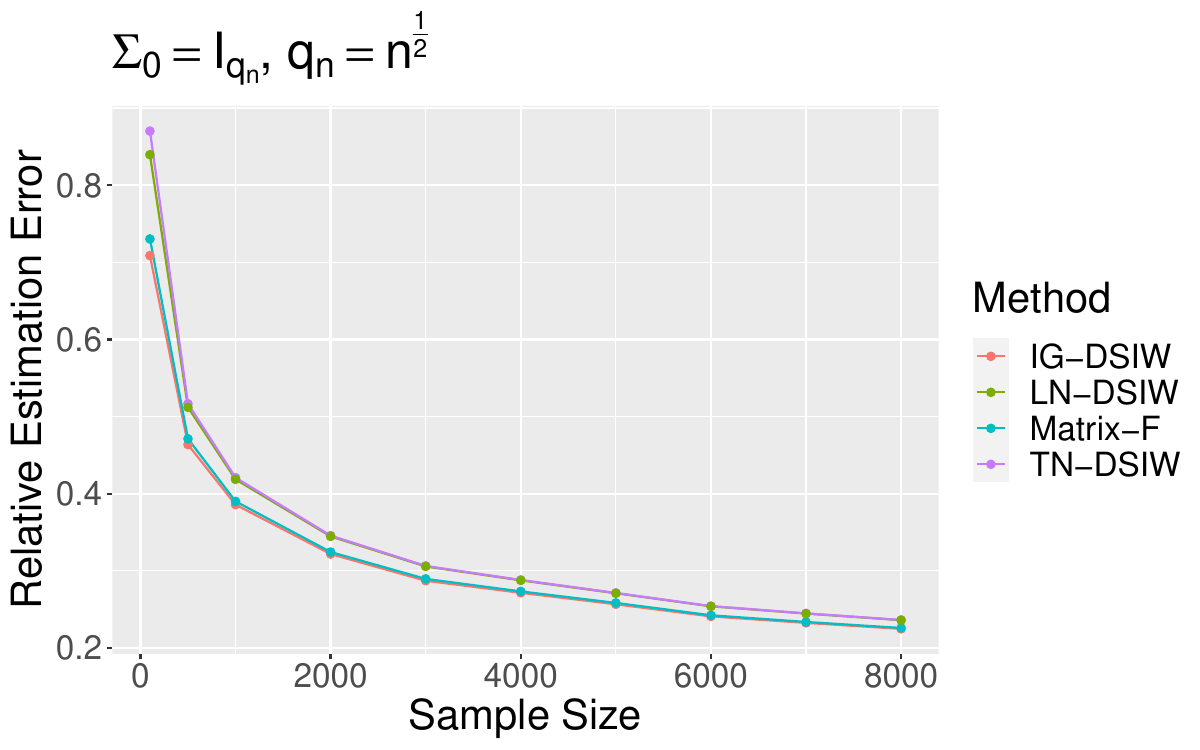}
  
\end{subfigure}
\begin{subfigure}{.4\textwidth}
  \centering
  \includegraphics[width=\linewidth]{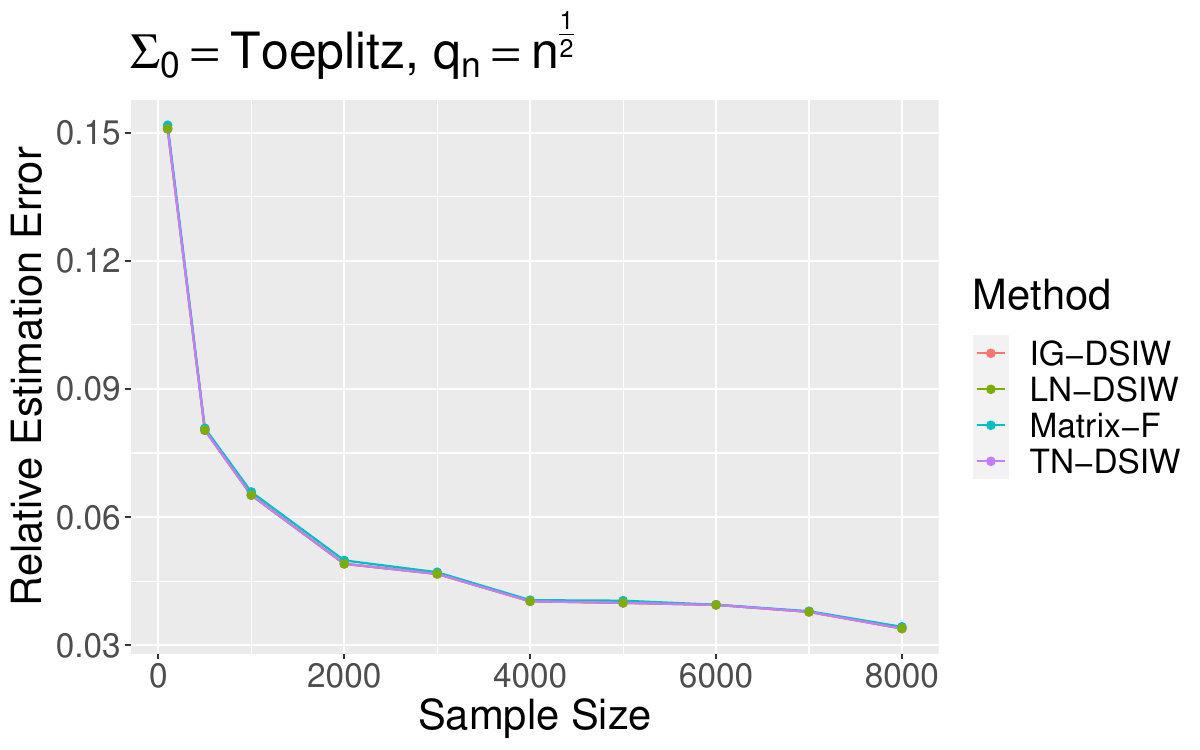}

\end{subfigure}%
\begin{subfigure}{.4\textwidth}
  \centering
  \includegraphics[width=\linewidth]{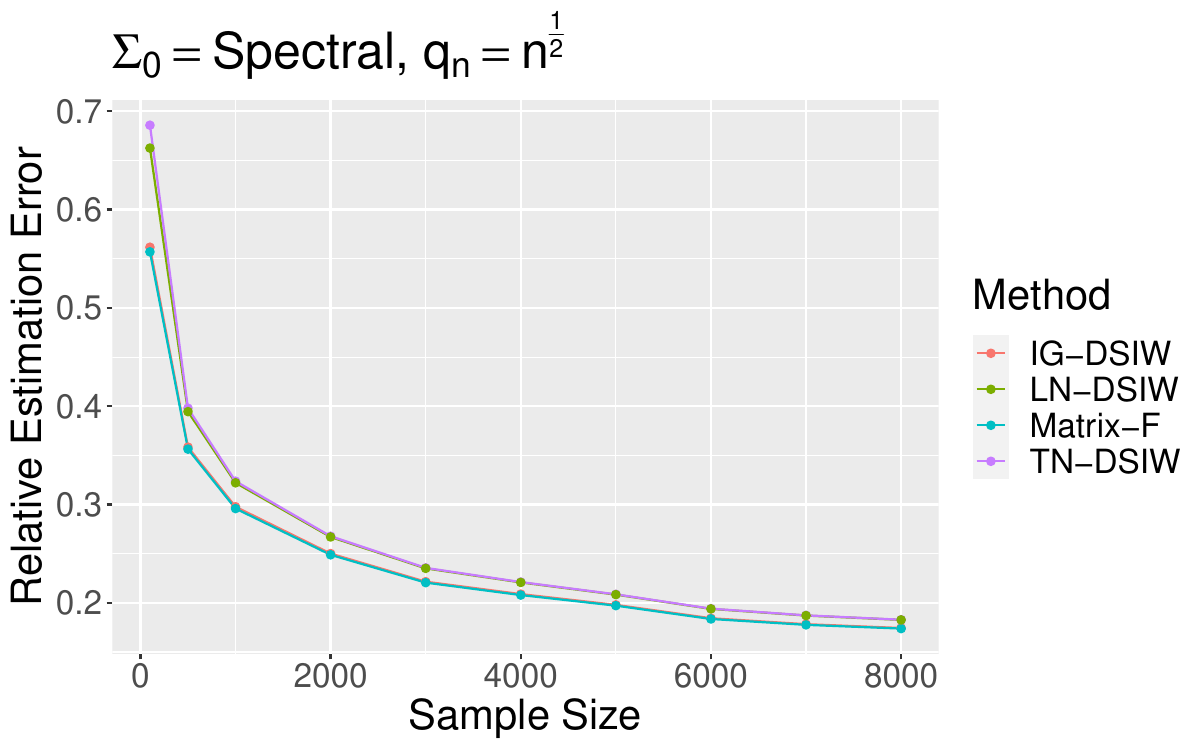}

\end{subfigure}
\begin{subfigure}{.4\textwidth}
  \centering
  \includegraphics[width=\linewidth]{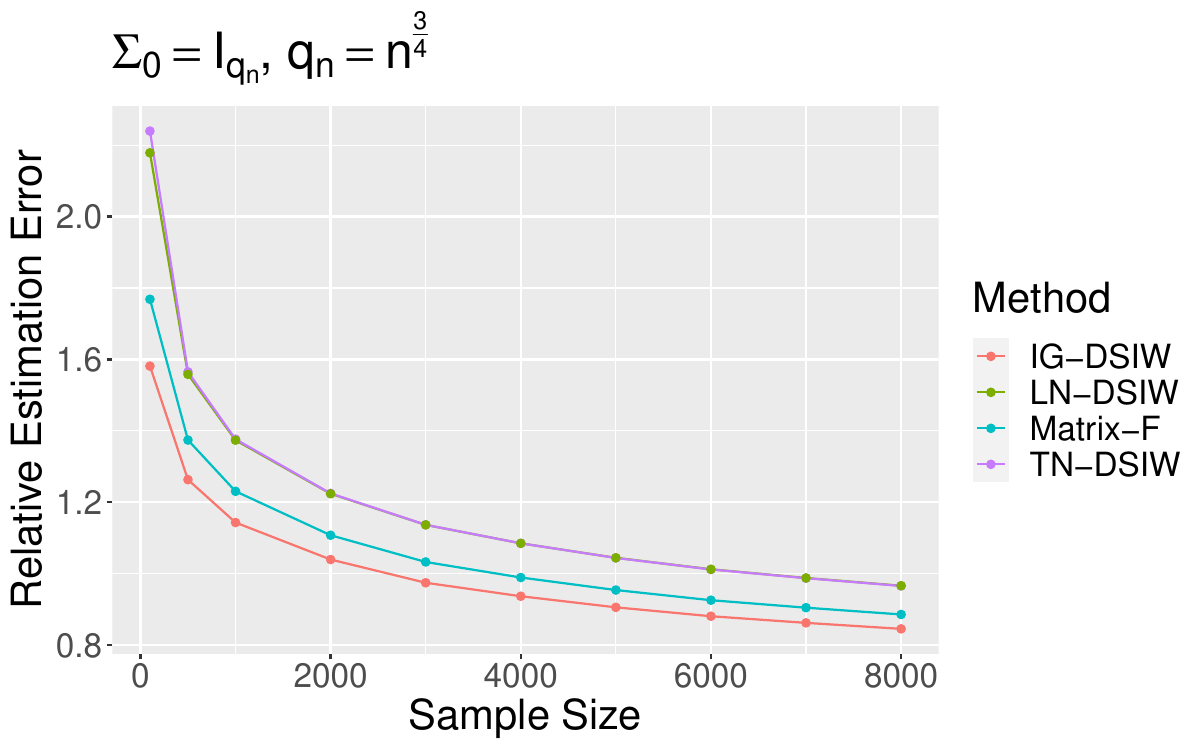}

\end{subfigure}%
\begin{subfigure}{.4\textwidth}
  \centering
  \includegraphics[width=\linewidth]{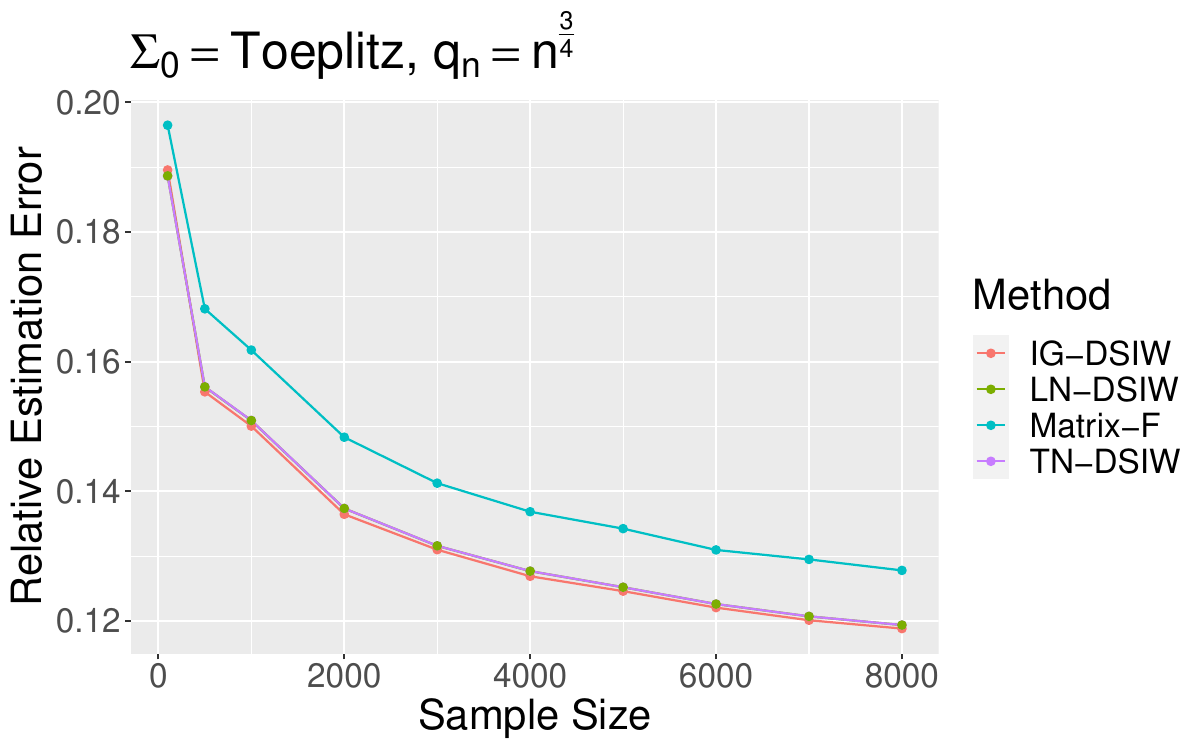}

\end{subfigure}
\begin{subfigure}{0.4\textwidth}
  \centering
  \includegraphics[width=\linewidth]{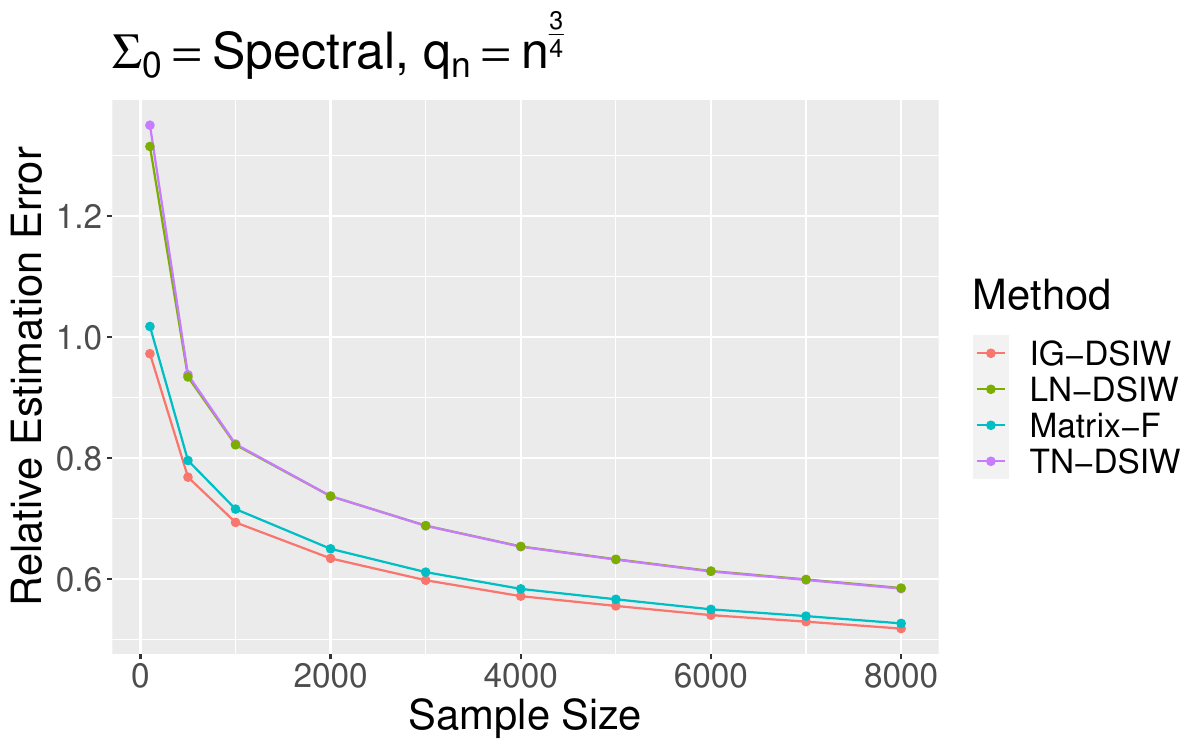}
 
\end{subfigure}
\caption{Plots of average relative estimation error for posterior mean $\hat{\uSigma}_n$ vs. sample size $n$ for nine different combinations of $\uSigma_0$ and $q_n$ (with $q_n/n \rightarrow 0$).}
\label{fig:fig2}
\end{figure}

\begin{figure}
\begin{subfigure}{0.3\textwidth}
  \centering
  \includegraphics[width=\linewidth]{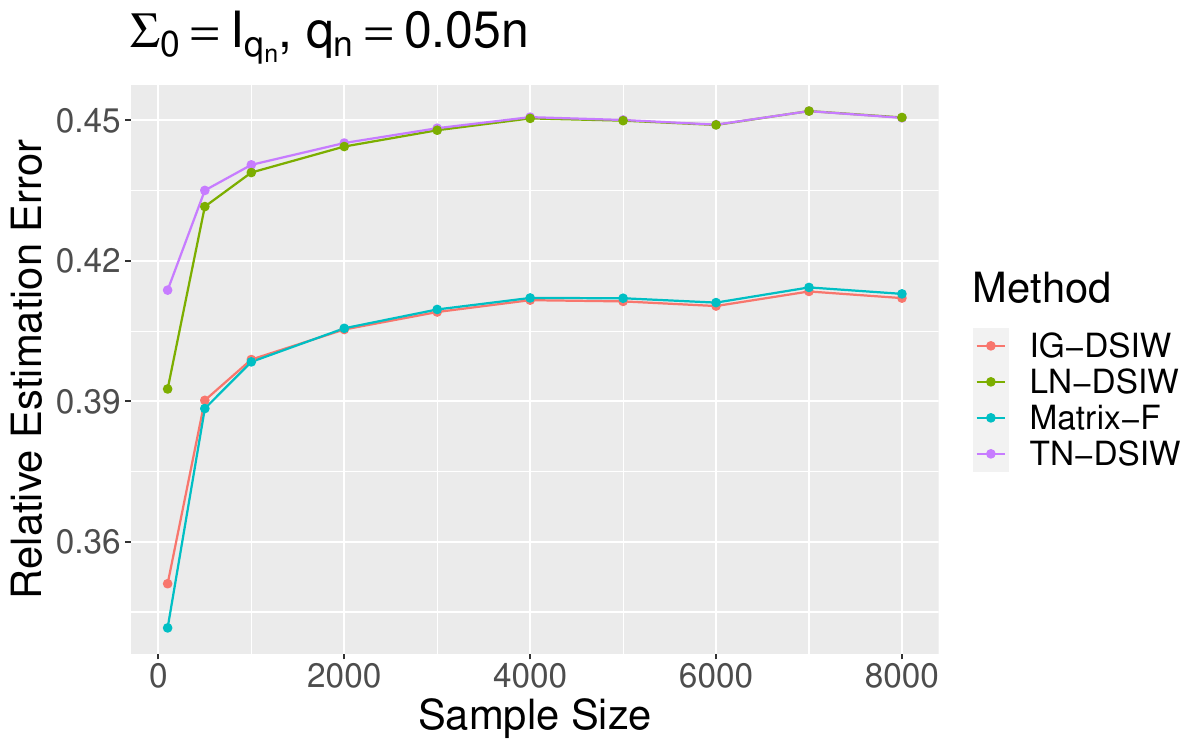}

\end{subfigure}%
\begin{subfigure}{0.3\textwidth}
  \centering
  \includegraphics[width=\linewidth]{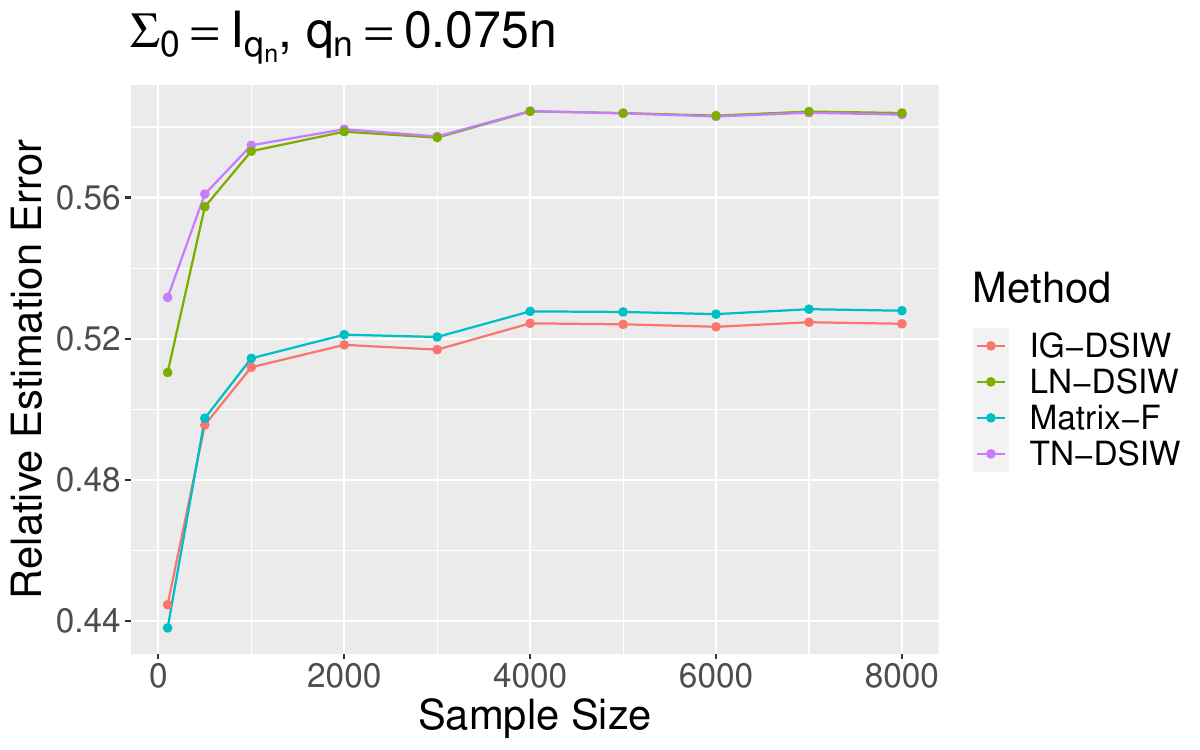}

\end{subfigure}
\begin{subfigure}{0.3\textwidth}
  \centering
  \includegraphics[width=\linewidth]{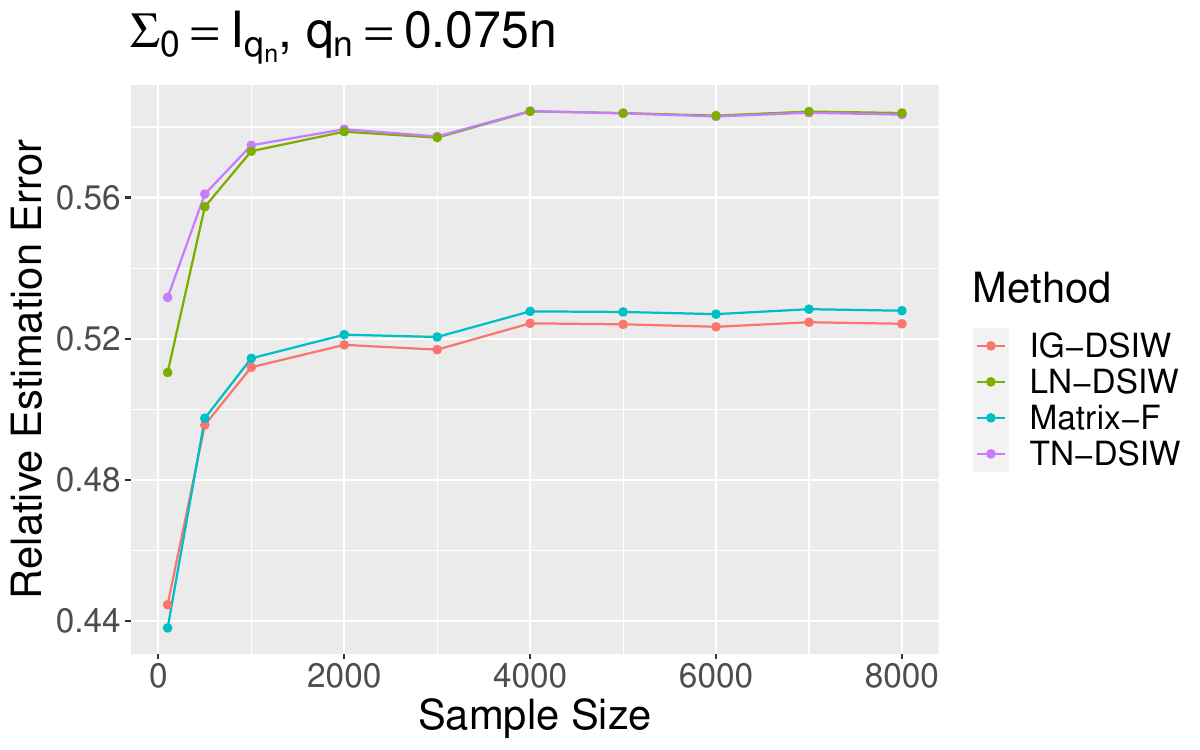}

\end{subfigure}%
\caption{Plots of average relative estimation error for posterior mean $\hat{\uSigma}_{n}$ vs. sample size $n$ with $\uSigma_0 = \uI_{q_n}$ and three different choices of $q_n$ (with $q_n$ increasing at same rate as $n$).}
\label{fig:fig3}
\end{figure}

\section{Simulation study}\label{sec9}
In this section, we conduct a simulation study to provide empirical evidence and demonstrate our posterior consistency and inconsistency results. Our primary focus lies on estimating the parameter $\uSigma$. Taking these factors into account and for the sake of simplicity, we adopt the i.i.d. covariance estimation setup outlined in Section \ref{sec7} solely for illustrative purposes.

\medskip

\noindent
{\it Simulation setup}. We use $10$ different values for the sample size $n$, ranging from $50$ to $8000$. We allow $q = q_n$ to increase with $n$ and consider three settings in particular: (a) $q_n=n^{\frac{1}{4}}$, (b) $q_n=n^{\frac{1}{2}}$, and (c) $q_n=n^{\frac{3}{4}}$. Three different choices of the true covariance matrix $\uSigma_0$ are considered, namely (a) $\uSigma=I_{q_n}$, 
(b) $\uSigma_0=((0.9^{|i-j|}))_{ij}$ (Toeplitz form), and (c) $\uSigma_0=\uU\uLambda\uU^T$ (spectral form). For the spectral form, $\uU^{q_n\times q_n}$ is a random orthogonal matrix and $\uLambda^{q_n\times q_n}$ is a diagonal matrix independent of $\uU$ whose diagonal entries are generated independently from Uniform$(1,2)$ distribution. For each $(n, q_n, \uSigma_0)$ combination, we generate $500$ datasets with $n$ i.i.d. observation from a multivariate Gaussian distribution with mean $0$ and covariance matrix $\uSigma_0$.

\medskip

\noindent{\it Evaluation strategy for posterior tail probabilities}. For each of the $500$ datasets we computed $\Pi_n(\norm{\uSigma-\uSigma_0} > M\sqrt{q_n/n}\mid\uYn)$ empirically for four different priors, namely IW-DSIW (\cite{huang}), LN-DSIW (\cite{zava}), TN-DSIW (\cite{zava}) and matrix-F (\cite{mulder}), with default hyperparameter choices. We used Gibbs sampling algorithms based on Lemmas \ref{Lemma_post1} and \ref{Lemma_post2}. In our Gibbs sampling chain of $10,000$ simulations, samples were thinned by keeping every 5th simulation after an initial burn-in period of $5000$ simulations, with the total number of samples kept being $1000$. Suppose $\{\uSigma^{(i)}\}_{i=1}^{1000}$ denotes the corresponding sequence of samples from the marginal posterior distribution of $\uSigma$. Then for an appropriately fixed constant $M > 0$, the quantity $$\frac{1}{1000}\sum_{i=1}^{1000}\uone_{\left\{\norm{\uSigma^{(i)}-\uSigma_0} > M\sqrt{q_n/n}\right\}},$$ 
is used as an empirical estimate of $\Pi_n(\norm{\uSigma-\uSigma_0} > M\sqrt{q_n/n}\mid\uYn)$. These quantities will be averaged over $500$ replicates to get a stable estimate of $\Pi_n(\norm{\uSigma-\uSigma_0} > M\sqrt{q_n/n}\mid\uYn)$. 

\medskip

\noindent{\it Evaluation strategy for posterior mean}. Again for each of the $500$ datasets, the posterior mean $\hat{\uSigma}_{n}$ is computed for four different priors, namely IW-DSIW (\cite{huang}), LN-DSIW (\cite{zava}), TN-DSIW (\cite{zava}) and matrix-F (\cite{mulder}), with default hyperparameter choices. The posterior means are computed using the aforementioned Gibbs sampling algorithms based on Lemmas \ref{Lemma_post1} and \ref{Lemma_post2}. The relative estimation errors of the posterior means, given by $\norm{\hat{\uSigma}_{n}-\uSigma_0}/\norm{\uSigma_0}$, are averaged over the $500$ replicates. 

\medskip

\noindent
{\it Results}. In the above setup, there are $9$ possible combinations of the pair ($\uSigma_0, q_n)$. For each of these combinations, Figure \ref{fig:fig1} provides a plot of the posterior tail probability vs. the sample size for each of the four prior choices mentioned above. For all $9$ plots, we see that the posterior tail probability approaches $0$ as the sample size increases. These results provide an empirical justification for Theorem \ref{th_siw} and \ref{th_matF} mentioned in Section \ref{sec4}.

Similarly, Figure \ref{fig:fig2} provides a plot of the relative estimation error vs. the sample size for each of the four prior choices mentioned above. For all 9 plots, we see that the relative 
estimation error decreases as the sample size increases, and the rate of decrease is sharper when $q_n/n\to 0$ at a faster rate. These results thereby provide empirical support for our consistency of posterior mean results in Section \ref{sec7}.  

For empirical validation of the inconsistency results in Section \ref{sec7}, we consider a similar simulation setup as above, but with $\uSigma_0=I_{q_n}$, and three choices of $q_n$ that grow at the same rate as $n$, namely $q_n=0.05n$, $q_n=0.075n$, and $q_n=0.1n$. The plots for the average relative estimation error for the posterior means corresponding to all four prior distributions vs. sample size are provided in Fig \ref{fig:fig3}. As the sample size increases, the relative estimation errors for all priors in all the plots stay away from $0$, and instead stabilize around a non-zero value as $n$ increases. These results thereby provide empirical support for our inconsistency 
results in Section \ref{sec7}.

\section{Proofs of the theorems from Sections \ref{sec4} and \ref{sec6}} \label{sec5}
We will now prove the results stated in Section \ref{sec4} and Section \ref{sec6}. The proofs make use of four technical lemmas which are stated below. The proofs of these lemmas are provided in the Supplementary Material (\cite{supp}). 

\begin{lemma}\label{subg}

 Let $\uP^{b\times b}$ be a projection matrix with rank $r(\leq b)$ and $\uZ^{b\times a}$ be a matrix whose rows are independent isotropic sub-Gaussian random vector with the sub-Gaussian parameter at most $\sigma_0$. Then $\uZ^{T}P\uZ$ can be written as $\uA^{T}\uA$, where $\uA^{r\times a}$ is a matrix whose all rows are independent isotropic sub-Gaussian random vector with the sub-Gaussian parameter at most $\tilde{c}\sigma_0$, where $\tilde{c}$ is a constant.
\end{lemma}

\begin{lemma}\label{subg2}
 Let $\uA$ be an $(n+q_n+\nu-1)\times q_n$ matrix whose entries are coming from an i.i.d $N(0,1)$ random variables with $q_n=o(n)$. Suppose, $c$ is an arbitrarily large constant. 
       With probability at least $1-2\;exp(- q_n/2)$ one has for sufficiently large n,
        $$\left(1+c \;\sqrt{\frac{q_n}{n}}\right)^{-\frac{1}{2}}\leq \frac{s_{min}(\uA)}{\sqrt{n+\nu-2}}\leq\frac{s_{max}(\uA)}{\sqrt{n+\nu-2}}\leq\left(1-c\; \sqrt{\frac{q_n}{n}}\right)^{-\frac{1}{2}}.$$
   \end{lemma}
   
\begin{lemma}\label{subg3}
 Let $\uA$ be an $(n-p_n)\times q_n$ matrix whose rows $A_i$ are independent sub-Gaussian isotropic random vectors with $\snorm{A_i}$ at most $\sigma_0$ with $p_n,q_n=o(n)$. Suppose, $c$ is an arbitrarily large constant and $c_{\sigma_{0}}$ is a constant that depends on $\sigma_0$. 
       With probability at least $1-2\;exp(- c_{\sigma_{0}}q_n/2)$ one has for sufficiently large n and for any sequence of positive real number $a_n\to 1$ as $n \to \infty$,
        $$\left(1-c\;a_n\sqrt{\frac{q_n}{n}}\right)^{\frac{1}{2}}\leq \frac{s_{min}(\uA)}{\sqrt{n-p_n}}\leq\frac{s_{max}(\uA)}{\sqrt{n-p_n}}\leq\left(1+c\;a_n \sqrt{\frac{q_n}{n}}\right)^{\frac{1}{2}}.$$
   \end{lemma}

\begin{lemma}\label{subg4}
 Suppose, in (\ref{model1}), $\BLS$ be the least square estimator of $\uB$, then under the true model $\mathbb{P}_0$ and sufficiently large constant $c$,
 \begin{align*}
     &\mathbb{P}_0\left(\norm{\BLS-\uB_0}\geq c \max\left(\sqrt{\frac{p_n}{n}},\;\sqrt{\frac{q_n}{n}}\right)\right)\leq2 \exp(-c_1\max{(\sqrt{p_n},\sqrt{q_n})}),\;\textit{and}\\
      &\mathbb{P}_0\left(\fnorm{\BLS-\uB_0}\geq c \sqrt{\frac{p_n q_n}{n}}\right)\leq2 \exp(-c_2p_nq_n),
 \end{align*}
 where $c_1$ and $c_2$ constants that depend on $c$.
 \end{lemma}

 \noindent
 With the lemmas in hand, we provide the proofs of Theorems \ref{th_siw} and 
 \ref{th_matF} in Sections \ref{proof_dsiw} and \ref{proof_matF} respectively. Additionally, the proofs of Theorems \ref{th_B1} and \ref{th_B2} can be found in Sections \ref{proof_B1} and \ref{proof_B2}.

\subsection{Proof of Theorem \ref{th_siw}} \label{proof_dsiw}

\noindent
To begin the proof first note that,
 \begin{align*}
    \Pi_n(\norm{\uSigma-\uSigma_0} > M\delta_n | \uYn)\leq &\;\Pi_n(\norm{\uSigma-\uSigma_{PM}} > M\delta_n/2 | \uYn)+\\
    &\;\Pi_n(\norm{\uSigma_{PM}-\uSigma_0} > M\delta_n/2 | \uYn),
 \end{align*}
where $\uSigma_{PM}:=E[\uSigma|\uDelta,\uYn]=\frac{c_{\nu}\uDelta}{n+\nu-2}+\frac{\uS_{Y}}{n+\nu-2}$, using posterior distributions in Lemma \ref{Lemma_post1}. Taking expectation on both sides, we get 
\begin{flalign}\label{th1:1}
    E_{\mathbb{P}_0}\left[\Pi_n(\norm{\uSigma-\uSigma_0} > M\delta_n | \uYn)\right]\leq &\;E_{\mathbb{P}_0}\left[\Pi_n(\norm{\uSigma-\uSigma_{PM}} > M\delta_n/2 | \uYn)\right]+\notag\\
    &E_{\mathbb{P}_0}\left[\Pi_n\left(\norm{\uSigma_{PM}-\uSigma_0} > M\delta_n/2 | \uYn\right)\right].
 \end{flalign}
To show the first term of the right side of the inequality goes to zero note that, if we define, $\uU^*=\uSigma_{PM}^{-1/2}\uSigma\uSigma_{PM}^{-1/2}$, then $\uU|\uDelta,\uYn\sim IW(n+\nu+q_n-1,(n+\nu-2)\uI_{q_n})$ from Lemma \ref{Lemma_post1}. Using the tower property of conditional expectations, we get 
\begin{flalign}\label{th1:2}
&E_{\mathbb{P}_0}\left[\Pi_n(\norm{\uSigma-\uSigma_{PM}} > M\delta_n/2 | \uYn)\right]\notag\\=&\; E_{\mathbb{P}_0}\left[\left[E_{\uDelta}\left[\Pi_n(\norm{\uSigma-\uSigma_{PM}} > M\delta_n/2 |\uDelta,\uYn)\right]\mid \uYn\right]\right]\notag\\=&\;
E_{\mathbb{P}_0}\left[E_{\uDelta}\left[\Pi_n\left(\norm{U-\uI_{q_n}} > \frac{M\delta_n}{2\sigma} |\uDelta,\uYn\right)\mid\uYn\right]\right]+E_{\mathbb{P}_0}\left[\Pi_n(\norm{\uSigma_{PM}} > \sigma | \uYn)\right].
\end{flalign}

\noindent
 Here $\sigma$ is a constant chosen such that $\sigma>2/k_{\sigma}$.  Now the second term of (\ref{th1:2}) will be smaller than 
 \begin{align*}
   E_{\mathbb{P}_0}\left[\Pi_n\left(\norm{\uSigma_{PM}-\uSigma_0} > \sigma/2 \mid \uY\right)\right]+\uone_{\left\{\norm{\uSigma_0}> \sigma/2\right\}}.   
 \end{align*}
 The indicator term is zero by the choice of $\sigma$ and note that,
\begin{align*}
    E_{\mathbb{P}_0}\left[\Pi_n(\norm{\uSigma_{PM}-\uSigma_0} > \sigma/2 | \uYn)\right]\leq E_{\mathbb{P}_0}\left[\Pi_n(\norm{\uSigma_{PM}-\uSigma_0} > M\delta_n/2 | \uYn)\right],
\end{align*}
for sufficiently large n since $q_n=o(n)$. Hence showing that the second term of (\ref{th1:1}) converges to zero will be sufficient to show that the second term of (\ref{th1:2}) converges to zero. 

For the first term of (\ref{th1:2}) note that,
$\norm{U-\uI_q}=\max\{|\lambda_{min}(\uU)-1|,|\lambda_{max}(\uU)-1|\}$. Also, by definition, conditional on $\uDelta$ and $\uYn,\;\uU^{-1}$  can be written as $\frac{\uA^{T}\uA}{n+\nu-2}$, when $\nu$ is an integer and $A$ is a $(\nu+q_n+n-1)\times q_n$ matrix where entries are coming from a i.i.d $N(0,1)$ distribution. For non-integer $\nu$ the same argument will follow by replacing $\nu$ by $\ceil{\nu}$ and $(\ceil{\nu}-1)$ and using the stochastic ordering property of the Wishart distribution with the increment and decrement in degrees of freedom. Now it follows that,
\begin{align}\label{th1:*1}
    &E_{\mathbb{P}_0}\left[E_{\uDelta}\left[\Pi_n\left(\norm{\uU-\uI_q} > \frac{M\delta_n}{2\sigma}|\uDelta,\uYn\right)|\uYn\right]\right]\notag\\\;\leq &\left\{1- E_{\mathbb{P}_0}\left[E_{\uDelta}\left[\Pi_n\left(\left(1+\frac{M\delta_n}{2\sigma}\right)^{-\frac{1}{2}}\leq \frac{s_{min}(\uA)}{\sqrt{n+\nu-2}}\leq\left(1-\frac{M\delta_n}{2\sigma}\right)^{-\frac{1}{2}}|\uDelta,\uYn\right)|\uYn\right]\right]\right\}+\notag\\&\;\left\{1- E_{\mathbb{P}_0}\left[E_{\uDelta}\left[\Pi_n\left(\left(1+\frac{M\delta_n}{2\sigma}\right)^{-\frac{1}{2}}\leq \frac{s_{max}(\uA)}{\sqrt{n+\nu-2}}\leq\left(1-\frac{M\delta_n}{2\sigma}\right)^{-\frac{1}{2}}|\uDelta,\uYn\right)|\uYn\right]\right]\right\}\notag\\
    \leq & \;4\exp\{-q_n/2\}    \to 0\;\textit{as}\;n \to \infty.
\end{align}
The final step follows from Lemma \ref{subg2} with $c=M/(2\sigma)$. 

We now analyze the second term of the right side of the inequality in (\ref{th1:1}). Note that
\begin{align}\label{th1:3}
    &E_{\mathbb{P}_0}\left[\Pi_n\left(\norm{\uSigma_{PM}-\uSigma_0} > M\delta_n/2 | \uYn\right)\right]\notag\\=&\;E_{\mathbb{P}_0}\left[\Pi_n\left(\norm{\frac{c_{\nu}\uDelta}{n+\nu-2}+\frac{\uS_{Y}}{n+\nu-2}-\uSigma_0} > M\delta_n/2 | \uYn\right)\right]\notag\\\leq&\;
    E_{\mathbb{P}_0}\left[\Pi_n\left(\norm{\frac{\uS_{Y}}{n+\nu-2}-\uSigma_0} > M\delta_n/4 | \uYn\right)\right]+E_{\mathbb{P}_0}\left[\Pi_n\left(\norm{\frac{c_{\nu}\uDelta}{n+\nu-2}} > M\delta_n/4 | \uYn\right)\right].
\end{align}
Recall that,
$\uS_Y=\uYn^{T}(\uI_n+\uXn\uX_{\lambda}^{-1}\uXn^{T})\uYn=\uW_n+\BLS^T(\lambda^{-1}\uI_{p_n}+(\uXn^{T}\uXn)^{-1})^{-1}\BLS$, using Woodbury's matrix identity where, $\uW_n=\uYn^{T}(\uI_n-\uXn(\uXn^{T}\uXn)^{-1}\uXn^{T})\uYn$. Thus the first term of the right side of the inequality in (\ref{th1:3}) can be written as
\begin{align}\label{th1:4}
    E_{\mathbb{P}_0}\left[\Pi_n\left(\norm{\frac{\uS_{Y}}{n+\nu-2}-\uSigma_0} > M\delta_n/4 | \uYn\right)\right]=\mathbb{P}_0\left(\norm{\frac{\uS_{Y}}{n+\nu-2}-\uSigma_0} > M\delta_n/4\right),
\end{align}
as $\uS_{Y}$ depends on data $\uYn$ only. Using the expression of $\uS_Y$ we can further reduce (\ref{th1:4}) as

\begin{align}\label{th1:5}
&\mathbb{P}_0\left(\norm{\frac{\uS_{Y}}{n+\nu-2}-\uSigma_0} > M\delta_n/4\right)  \notag\\\leq & \;\mathbb{P}_0\left(\norm{\frac{\uWn}{n+\nu-2}-\uSigma_0} > M\delta_n/8\right)+\notag\\&\;
\mathbb{P}_0\left(\sqrt{\norm{(\lambda^{-1}\uI_{p_n}+(\uXn^{T}\uXn)^{-1})^{-1}}}\times\norm{\BLS}  >\sqrt{ \frac{{M\delta_n(n+\nu-2)}}{8}}\right).
\end{align}
First part of (\ref{th1:5}) can be manipulated as
\begin{align}\label{th1:6}
    \mathbb{P}_0\left(\norm{\frac{\uWn}{n+\nu-2}-\uSigma_0} > M\delta_n/8\right)\leq &\;\mathbb{P}_0\left(\norm{\frac{\uWn}{n-p_n}-\uSigma_0} >  M a_n\delta_n/16\right)+\notag\\
    &\; \uone_{\left\{\frac{\nu+p_n-2}{n+\nu-2}\norm{\uSigma_0}>M\delta_n/16\right\}},
\end{align}
where $a_n=(n+\nu-2)/(n-q_n)$. Using Assumptions  \ref{as2} and \ref{as4} for sufficiently large $n$, the second term of (\ref{th1:6}) is $0$. Also
\begin{align*}
    \uWn=\uYn^{T}(\uI_n-\uXn(\uXn^{T}\uXn)^{-1}\uXn^{T})\uYn=&\;\uE_n^{T}(\uI_n-\uXn(\uXn^{T}\uXn)^{-1}\uXn^{T})\uE_n\\=&\;{\uSigma_0}^{1/2}{\uE_{n}^{*}}^T\uXn(\uXn^{T}\uXn)^{-1}\uXn^T\uE_{n}^{*}{\uSigma_0}^{1/2},
\end{align*}
where $\uE_{n}^{*}=\uE_{n}{\uSigma_0}^{-1/2}$. Then each row of $\uE_{n}^{*}$ are independent isotropic sub-Gaussian random variable with the sub-Gaussian norm at most $\sigma_0$. Since $(\uI_n-\uXn(\uXn^{T}\uXn)^{-1}\uXn^{T})$ is the projection matrix of rank $(n-p_n)$ we can write ${\uE_{n}^{*}}^T\uXn(\uXn^{T}\uXn)^{-1}\uXn^T\uE_{n}^{*}=\uA_1^{T}\uA_1,$
where $\uA_1^{n-p_n\times q_n}$ is a matrix whose all rows are independent isotropic sub-Gaussian random vector with the sub-Gaussian parameter at most $\tilde{c}\sigma_0$, where $\tilde{c}$ is a constant using Lemma \ref{subg} with $a=q_n,\;b=n$ and $r=n-p_n$. Keeping Assumption \ref{as2} in mind, the first term of (\ref{th1:6}) can be written as 
\begin{align*}
    & \mathbb{P}_0\left(\norm{\frac{\uWn}{n-p_n}-\uSigma_0}>Ma_n \delta_n/16\right) \leq \mathbb{P}_0\left(\norm{\frac{\uA_1^{T}\uA_1}{n-p_n}-\uI_{q_n}}>k_{\sigma} a_n M\delta_n/16\right)\\\leq &\;1- \mathbb{P}_0\left(\left(1-Ma_n \delta_n/16\right)^{\frac{1}{2}}\leq \frac{s_{min}(\uA_1)}{\sqrt{n-p_n}}\leq\left(1+Ma_n \delta_n/16\right)^{\frac{1}{2}}\right)+\\&\;1-\mathbb{P}_0\left(\left(1-Ma_n \delta_n/16\right)^{\frac{1}{2}}\leq \frac{s_{max}(\uA_1)}{\sqrt{n-p_n}}\leq\left(1+Ma_n \delta_n/16\right)^{\frac{1}{2}}\right)\leq4\;exp(-\frac{c_{\sigma_{0}}q_n}2)\\&\;
    \to 0\; \textit{as}\; n \to \infty.
\end{align*}
The last step follows from Lemma \ref{subg3} with $c=M/16$. 

It now remains to be shown that the second term of (\ref{th1:3}) and the second term of (\ref{th1:5}) 
converge to zero. Since 
$
\norm{(\lambda^{-1}\uI_{p_n}+ (\uXn^{T}\uXn)^{-1})^{-1}}\leq \norm{\lambda\uI_{p_n}}=\lambda,
$
using Assumption \ref{as5}, the second term of (\ref{th1:5}) can be written as

\begin{align}\label{th1:7}
    &\mathbb{P}_0\left(\sqrt{\norm{(\lambda^{-1}\uI_{p_n}+(\uXn^{T}\uXn)^{-1})^{-1}}}\times\notag\norm{\BLS}  >\sqrt{ \frac{{M\delta_n(n+\nu-2)}}{8}}\right)\notag\\
    \leq&\;\mathbb{P}_0\left(\norm{\BLS-\uB_0}  >\sqrt{ \frac{{M\lambda^2\delta_n(n+\nu-2)}}{32}}\right)+ \notag\\
    &\;\uone_{\left\{\sqrt{\norm{(\lambda^{-1}\uI_{p_n}+(\uXn^{T}\uXn)^{-1})^{-1}}}\norm{\uB_0}>\sqrt{ \frac{{M\delta_n(n+\nu-2)}}{32}}\right\}}.
\end{align}
The first term of (\ref{th1:7}) will go to $0$ as $n \to \infty$ because of the fact $\delta_n(n+\nu-2) \to \infty$ as $n \to \infty$ using Lemma \ref{subg4}. For the next part observe that
$\norm{(\lambda^{-1}\uI_{p_n}+(\uXn^{T}\uXn)^{-1})^{-1}}\leq (\lambda^{-1}+1/\norm{(\uXn^{T}\uXn)})^{-1}\leq((\lambda^{-1}+k_x/n)^{-1},$
 using Assumption \ref{as5}. Thus we can manipulate the condition within the indicator function and write it as
 \begin{align*}
     \uone_{\left\{\sqrt{\norm{(\lambda^{-1}\uI_{p_n}+(\uXn^{T}\uXn)^{-1})^{-1}}}\norm{\uB_0}>\sqrt{ \frac{{M\delta_n(n+\nu-2)}}{32}}\right\}} \leq  \uone_{\left\{\lambda^{-1}\leq\frac{32\norm{{\uB_0}}^2}{M\delta_n(n+\nu-2)}-\frac{k_x}{n}\right\}}.
 \end{align*}
The last indicator term can be made $0$ using Assumption \ref{as6} for an appropriate value of $\lambda_0$. 

Finally, we have to deal with the second term of (\ref{th1:3}). Define $\omega_n:=M\delta_n(n+\nu-2)/(4c_\nu)$, then it follows that
\begin{align}\label{th1:8}
    \Pi_n(\norm{\uDelta} > \omega_n | \uYn)=  \Pi_n(\max_{1\leq i\leq q_n}{\delta_{i}} > \omega_n | \uYn)\leq \sum_{i=1}^{q_n}\Pi_n({\delta_{i}} > \omega_n | \uYn).
\end{align}
If any of the $\delta_i$'s has bounded support then for sufficient large $n,\;\Pi_n({\delta_{i}} > \omega_n | \uYn)$ will be $0$ and there is nothing to show. Hence without loss of generality, we will assume $\delta_i$'s has a positive density over the whole positive real line. Using the form of posterior distributions from Lemma \ref{Lemma_post1} it follows that,

\begin{align}\label{th1:9}
    \Pi_n({\delta_{1}} > \omega_n | \uYn)=&\;\frac{\int_{\omega_n}^\infty\int_{0}^\infty\dots\int_{0}^\infty\left|\uS_Y+c_{\nu}\uDelta\right|^{-\frac{n+\nu+q_n-1}{2}}\prod_{i=1}^{q_n} \delta_i^{\frac{\nu+q_n-1}{2}}\pi_i(\delta_i)  \,\prod_{i=1}^{q_n}d\delta_{i}}{\int_{0}^\infty\dots\int_{0}^\infty\left|\uS_Y+c_{\nu}\uDelta\right|^{-\frac{n+\nu+q_n-1}{2}}\prod_{i=1}^{q_n} \delta_i^{\frac{\nu+q_n-1}{2}}\pi_i(\delta_i)  \,\prod_{i=1}^{q_n}d\delta_{i}}\notag\\
    \leq&\;\frac{\int_{\omega_n}^\infty\int_{0}^\infty\dots\int_{0}^\infty\left|\uS_Y+c_{\nu}\uDelta\right|^{-\frac{n+\nu+q_n-1}{2}}\prod_{i=1}^{q_n}
    \delta_i^{\frac{\nu+q_n-1}{2}}\pi_i(\delta_i)  \,\prod_{i=1}^{q_n}d\delta_{i}}{\int_0^{{\omega_n}}\int_{0}^\infty\dots\int_{0}^\infty\left|\uS_Y+c_{\nu}\uDelta\right|^{-\frac{n+\nu+q_n-1}{2}}\prod_{i=1}^{q_n} \delta_i^{\frac{\nu+q_n-1}{2}}\pi_i(\delta_i)  \,\prod_{i=1}^{q_n}d\delta_{i}}.
\end{align}

From this step, we will assume $c_\nu=1$ to keep our steps simple since degrees of freedom are a fixed value in our case. Write  $\left|\uS_Y+\uDelta\right|=\left|diag(\delta_1,\dots\delta_p)+
  \begin{bmatrix}
    s_{y1} & u^{T}  \\
    u & C^{*} 
  \end{bmatrix}
\right|$, where $s_{y1}$ is the first diagonal entry of $\uS_Y$ and $u$ and $C^{*}$ be the corresponding entries. By using the Schur complement determinant formula note that $\left|\uS_Y+\uDelta\right|=(s_{y1}+\delta_1)\left|C-\frac{uu^{T}}{s_{y1}+\delta_1}\right|$, where $C=C^{*}+diag(\delta_2,\dots\delta_p)$. So, $C$ and $u$ is completely free from $\delta_1$. Also $\left|C-\frac{uu^{T}}{s_{y1}+\delta_1}\right|$ is a positive definite matrix and is an increasing function of $\delta_1$. Keeping these and, Assumption \ref{as3} in mind write (\ref{th1:9}) as
\begin{align}\label{th1:10}
    &\Pi_n({\delta_{1}} > \omega_n | \uYn)\notag\\
    \leq&\;\frac{\int_{\omega_n}^\infty\int_{0}^\infty\dots\int_{0}^\infty({s_{y1}+\delta_1})^{-\frac{n+\nu+q_n-1}{2}}\left|C-\frac{uu^{T}}{s_{y1}+\delta_1}\right|^{-\frac{n+\nu+q_n-1}{2}}\prod_{i=1}^{q_n} \delta_i^{\frac{\nu+q_n-1}{2}}\pi_i(\delta_i)  \,\prod_{i=1}^{q_n}d\delta_{i}}{\int_0^{{\omega_n}}\int_{0}^\infty\dots\int_{0}^\infty({s_{y1}+\delta_1})^{-\frac{n+\nu+q_n-1}{2}}\left|C-\frac{uu^{T}}{s_{y1}+\delta_1}\right|^{-\frac{n+\nu+q_n-1}{2}}\prod_{i=1}^{q_n} \delta_i^{\frac{\nu+q_n-1}{2}}\pi_i(\delta_i)  \,\prod_{i=1}^{q_n}d\delta_{i}}\notag\\
    \leq&\; \frac{\int_{\omega_n}^\infty\int_{0}^\infty\dots\int_{0}^\infty({s_{y1}+\delta_1})^{-\frac{n+\nu+q_n-1}{2}}\left|C-\frac{uu^{T}}{s_{y1}+\omega_n}\right|^{-\frac{n+\nu+q_n-1}{2}}\prod_{i=1}^{q_n} \delta_i^{\frac{\nu+q_n-1}{2}}\pi_i(\omega_n)  \,\prod_{i=1}^{q_n}d\delta_{i}}{\int_k^{{\omega_n}}\int_{0}^\infty\dots\int_{0}^\infty({s_{y1}+\delta_1})^{-\frac{n+\nu+q_n-1}{2}}\left|C-\frac{uu^{T}}{s_{y1}+\omega_n}\right|^{-\frac{n+\nu+q_n-1}{2}}\prod_{i=1}^{q_n} \delta_i^{\frac{\nu+q_n-1}{2}}\pi_i(\omega_n)  \,\prod_{i=1}^{q_n}d\delta_{i}}\notag\\
    \leq&\;\frac{\int_{\omega_n}^\infty({s_{y1}+\delta_1})^{-\frac{n+\nu+q_n-1}{2}}\delta_1^{\frac{\nu+q_n-1}{2}}  \,d\delta_{1}}{\int_k^{\omega_n}({s_{y1}+\delta_1})^{-\frac{n+\nu+q_n-1}{2}}\delta_1^{\frac{\nu+q_n-1}{2}}} = \frac{\int_{\tilde{\omega}_n}^1 u_1^{\frac{\nu+q_n-1}{2}}(1-u_1)^{\frac{n}{2}-2}du_1}{\int_{k_n}^{\tilde{\omega}_n} u_1^{\frac{\nu+q_n-1}{2}}(1-u_1)^{\frac{n}{2}-2}du_1},
\end{align}
using the substitution $u_1=\delta_1/(\delta_1+s_{y1})$, where $\tilde{\omega}_n=\omega_n/(\omega_n+s_{y1})$ and $k_n=k/(k+s_{y1})$. Note that $s_{y1}$ must be in between $\lambda_{min}(\uS_Y)$ and $\lambda_{max}(\uS_Y)$. Thus using (\ref{th1:10}), (\ref{th1:8}) reduces to

\begin{align}\label{th1:11}
     \Pi_n(\norm{\uDelta} > \omega_n | \uYn)\leq \sum_{i=1}^{q_n}\Pi_n({\delta_{i}} > \omega_n | \uYn) \leq q_n \frac{P(U_n>\omega_n^*)}{1-P(U_n>\omega_n^*)-P(U_n<k_n^*)},
\end{align}
where $U_n \sim Beta(\frac{\nu+q_n-1}{2},\frac{n-2}{2})$ and $\omega_n^*=\omega_n/(\omega_n+\lambda_{max}(\uS_Y))$ and $k_n^*=k/(k+\lambda_{min}(\uS_Y))$. Consider the set,
$A_n=\{\norm{\uS_Y/n-\uSigma_0}<M\delta_n\}$ then we know $\mathbb{P}_0(A_n)\to1$ as $n \to \infty$ using similar derivation as in (\ref{th1:5}). Observe in the set $A_n$, $\lambda_{max}(\uS_Y/n)\leq M\delta_n + \norm{\uSigma_0}\leq M\delta_n+(1/k_{\sigma})$, using Assumption \ref{as2}. Also in a set $B_n=\{\lambda_{min}({W_n}/n)\geq(1-\eta)\lambda_{min}(\uSigma_0)\}$ we have $\lambda_{min}(\uS_Y/n)\geq\lambda_{min}({W_n}/n)\geq(1-\eta)\lambda_{min}(\uSigma_0)\geq(1-\eta)k_{\sigma}$ for sufficiently large n, using Weyl's inequality and Assumption \ref{as2} where $\eta$ is fixed constant in $(0,1)$. From the remark mentioned in the supplementary material (\cite{supp}) (after the proof of Lemma \ref{subg3}) we know $\mathbb{P}_0(B_n)\to1$ as $n \to \infty$. Defining $\bar{\omega}_n:=\frac{(\omega_n/n)}{((\omega_n/n)+M\delta_n+(1/k_{\sigma}))}$ and $\bar{k}_n:=\frac{k}{k+n(1-\eta)k_{\sigma}}$. Thus for sufficiently large $n$ within the set $A_n^*=(A_n\cap B_n)$, (\ref{th1:11}) reduces to
\begin{align}\label{th1:12}
     \Pi_n(\norm{\uDelta} > \omega_n | \uYn) \leq q_n\; \frac{P(U_n>\bar{\omega}_n)}{1-P(U_n>\bar{\omega}_n)-P(U_n<\bar{k}_n)}.
\end{align}
Since $E(U_n)=\frac{\nu+q_n+1}{\nu+q_n+n-1}$ and $Var(U_n)=\frac{(\nu+q_n+1)(n-2)}{(\nu+q_n+n-1)(\nu+q_n+n+1)^2}\leq\frac{\nu+q_n+1}{n^2}$ we can get $(\bar{\omega}_n-E(U))\geq \sqrt{c_1 \frac{\nu+q_n+1}{n}}>0$ for a constant $c_1$ and sufficiently large $n$. Thus using Chebychev's inequality
\begin{align}\label{th1:13}
    P(U_n>\bar{\omega}_n)=&\;P(U_n-E(U_n)>\bar{\omega}_n-E(U_n))\notag\\
    \leq &\;\frac{Var(U_n)}{(\bar{\omega}_n-E(U_n))^2}\leq\frac{1}{nc_1}.
\end{align}
Also $(\bar{k}_n-E(U))\leq -\sqrt{\frac{c_2}{n}}<0$ for a constant $c_2$ and sufficiently large $n$. Again by using Chebychev's inequality as in (\ref{th1:13}), one can show
\begin{align}\label{th1:14}
    P(U_n<k_n)\leq\frac{\nu+q_n+1}{nc_2}.
\end{align}
Thus using (\ref{th1:13}) and (\ref{th1:14}) within the set $A_n^*$
\begin{align}\label{th1:15}
   \Pi_n(\norm{\uDelta} > \omega_n | \uYn) \leq \frac{\frac{q_n}{nc_1}}{1-\frac{1}{nc_1}-\frac{\nu+q_n+1}{nc_2}}\to 0\;\textit{as}\; n \to \infty,
\end{align}
since $q_n=o(n)$ by Assumption \ref{as1}. Now we conclude the whole proof since $\mathbb{P}_0(A_n^{*c})\to 0$ as $n \to \infty$ implies  $E_{\mathbb{P}_0}\left[\Pi_n(\norm{\uDelta} > \omega_n | \uYn)\right] \to 0$ as $n \to \infty$.

\subsection{Proof of Theorem \ref{th_matF}} \label{proof_matF}

As discussed in Section \ref{sec3} the main difference between DSIW prior and Matrix-$F$ prior is the scale parameter for the base Inverse Wishart distribution i.e. $\bar{\uDelta}$ is now a general positive definite matrix. Hence the initial steps for this proof will be similar to the proof of Theorem \ref{th_siw} for the DSIW prior. We will only discuss the proof more rigorously when the posterior distribution of $\bar{\uDelta}\mid\uYn$ comes into play. To begin the proof first note that,
 \begin{align*}
    \Pi_n(\norm{\uSigma-\uSigma_0} > M\delta_n | \uYn)\leq &\;\Pi_n(\norm{\uSigma-\uSigma_{PM}} > M\delta_n/2 | \uYn)+\\
    &\;\Pi_n(\norm{\uSigma_{PM}-\uSigma_0} > M\delta_n/2 | \uYn),
 \end{align*}
where $\uSigma_{PM}:=E[\uSigma|\bar{\uDelta},\uYn]=\frac{\bar{\uDelta}}{n+\nu-2}+\frac{\uS_{Y}}{n+\nu-2}$, using posterior distributions in Lemma \ref{Lemma_post2}. Taking expectations on both sides, we get 
\begin{flalign}\label{th2:*1}
    E_{\mathbb{P}_0}\left[\Pi_n(\norm{\uSigma-\uSigma_0} > M\delta_n | \uYn)\right]\leq &\;E_{\mathbb{P}_0}\left[\Pi_n(\norm{\uSigma-\uSigma_{PM}} > M\delta_n/2 | \uYn)\right]+\notag\\
    &E_{\mathbb{P}_0}\left[\Pi_n\left(\norm{\uSigma_{PM}-\uSigma_0} > M\delta_n/2 | \uYn\right)\right].
 \end{flalign}
To show the first term of the right side of the inequality goes to zero note that, if we define, $\uU^*=\uSigma_{PM}^{-1/2}\uSigma\uSigma_{PM}^{-1/2}$, then $\uU^*|\bar{\uDelta},\uYn\sim IW(n+\nu+q_n-1,(n+\nu-2)\uI_{q_n})$ from Lemma \ref{Lemma_post2}. Using the similar steps as in (\ref{th1:2}) of Theorem \ref{th_siw}, we get 
\begin{flalign}\label{th2:*2}
&E_{\mathbb{P}_0}\left[\Pi_n(\norm{\uSigma-\uSigma_{PM}} > M\delta_n/2 | \uYn)\right]\notag\\
    \leq&\; E_{\mathbb{P}_0}\left[E_{\bar{\uDelta}}\left[\Pi_n\left(\norm{U^*-\uI_{q_n}} > \frac{M\delta_n}{2\sigma} |\bar{\uDelta},\uYn\right)\mid\uYn\right]\right]+E_{\mathbb{P}_0}\left[\Pi_n(\norm{\uSigma_{PM}} > \sigma | \uYn)\right].
\end{flalign}

\noindent
 Here $\sigma$ is a constant chosen such that $\sigma>2/k_{\sigma}$. Now we can claim the first term of \ref{th2:*2} converges to zero by using similar steps as in (\ref{th1:*1}) of Theorem \ref{th_siw}. Now the second term of (\ref{th2:*2}) will be smaller than $E_{\mathbb{P}_0}\left[\Pi_n(\norm{\uSigma_{PM}-\uSigma_0} > \sigma/2 | \uYn)\right]+\uone_{\left\{\norm{\uSigma_0}> \sigma/2\right\}}$. The indicator term is zero by the choice of $\sigma$ and note that,
\begin{align*}
    E_{\mathbb{P}_0}\left[\Pi_n(\norm{\uSigma_{PM}-\uSigma_0} > \sigma/2 | \uYn)\right]\leq E_{\mathbb{P}_0}\left[\Pi_n(\norm{\uSigma_{PM}-\uSigma_0} > M\delta_n/2 | \uYn)\right],
\end{align*}
for sufficiently large n since $q_n=o(n)$. Hence showing that the second term of (\ref{th2:*1}) converges to zero will be sufficient to show that the second term of (\ref{th2:*2}) converges to zero.

Again for the second term of (\ref{th2:*1}) using (\ref{th1:3}) of Theorem \ref{th_siw} we can write,
\begin{align}\label{th2:*3}
    &E_{\mathbb{P}_0}\left[\Pi_n\left(\norm{\uSigma_{PM}-\uSigma_0} > M\delta_n/2 | \uYn\right)\right]\notag\\\leq&\;
    \mathbb{P}_0\left(\norm{\frac{\uS_{Y}}{n+\nu-2}-\uSigma_0} > M\delta_n/4\right)+E_{\mathbb{P}_0}\left[\Pi_n\left(\norm{\frac{\bar{\uDelta}}{n+\nu-2}} > M\delta_n/4 | \uYn\right)\right].
\end{align}
Now it is easy to see  $\mathbb{P}_0\left(\norm{\frac{\uS_{Y}}{n+\nu-2}-\uSigma_0} > M\delta_n/4\right)\to 0$ as $n \to \infty$ using similar steps as in (\ref{th1:5}-\ref{th1:7}) of Theorem \ref{th_siw}. Now, all we need to show the second term of (\ref{th2:*3}) converges to zero as $n \to \infty$. 

To show $E_{\mathbb{P}_0}\left[\Pi_n(\norm{\bar{\uDelta}} > 2\omega_n | \uYn)\right]\to 0$ as $n \to \infty$, where $\omega_n:=M\delta_n(n+\nu-2)/(8)$, we must keep in mind that our $\bar{\uDelta}$ is now a general positive definite matrix (instead of a diagonal matrix $\uDelta$ in Theorem \ref{th_siw}). So the problem is more complicated now and we need different mathematical notions to deal with it. 

Let $\mathcal{N}_{\frac{1}{4}}$ be a $\frac{1}{4}$-net of $\mathcal{S}^{q_n-1}$ and say $\mathcal{N}(\mathcal{S}^{q_n-1}, {\frac{1}{4}})$ is the minimal cardinality of $\mathcal{N}_{\frac{1}{4}}$. From \cite[Lemma 5.2]{vershynin} we know $\mathcal{N}(\mathcal{S}^{q_n-1}, {\frac{1}{4}})\leq9^{q_n}$. Also from  From \cite[Lemma 5.3]{vershynin} it follows that $\norm{\bar{\uDelta}}\leq 2\underset{{u\in \mathcal{N}_{\frac{1}{4}}}}{\sup}u^T\bar{\uDelta}u$. Thus
\begin{align}\label{th2:2}
    \Pi_n(\norm{\bar{\uDelta}} > 2\omega_n | \uYn) \leq &\;\Pi_n(\underset{{u\in \mathcal{N}_{\frac{1}{4}}}}{\sup}u^T\bar{\uDelta}u > \omega_n | \uYn)
    \leq 9^{q_n} \underset{\norm{u}= 1}{\sup}\Pi_n(u^T\bar{\uDelta}u > \omega_n|\uYn),
\end{align}
where last inequality follows by taking union over all vectors $u \in \mathcal{N}_{\frac{1}{4}}$. Given any $u$ with $\norm{u}=1$, suppose $\uH_u^{q_n\times q_n}$ be a orthogonal matrix whose first column is $u$. Then Jacobian of transformation from $\bar{\uDelta}\to\tilde{\uDelta}=\uH_u^T\bar{\uDelta}\uH$ is $|\uH_u|^{q_n+1}=1$ using \cite[Proposition 5.11]{eaton}. Then if $\tilde{\Delta}_{11}$ is the first diagonal element of $\tilde{\uDelta}$ using the form of marginal posterior distribution of $\bar{\uDelta}$ from Lemma \ref{Lemma_post2} it follows,
\begin{align}\label{th2:3}
    &\Pi_n(u^T\bar{\uDelta}u > \omega_n| \uYn)=\Pi_n(\tilde{\Delta}_{11}> \omega_n| \uYn)\notag\\
       =&\; \frac{\int_{\mathbb{P}^{+}_{q_n}}\uone(\tilde{\Delta}_{11}> \omega_n) \left|\uS_{Yu}+\tilde{\uDelta}\right|^{-\frac{n+\nu+q_n-1}{2}}|\tilde{\uDelta}|^{\frac{\nu+\nu_q^*-2}{2}}\exp(-\tr(\uPsi^{-1}_u\tilde{\uDelta})/2)d\tilde{\uDelta}}
       {\int_{\mathbb{P}^{+}_{q_n}} \left|\uS_{Yu}+\tilde{\uDelta}\right|^{-\frac{n+\nu+q_n-1}{2}}|\tilde{\uDelta}|^{\frac{\nu+\nu_q^*-2}{2}}\exp(-\tr(\uPsi^{-1}_u\tilde{\uDelta})/2)d\tilde{\uDelta}},
\end{align}
where $\uS_{Yu}=\uH_u^T\uS_{Y}\uH_u$ and $\uPsi_u=\uH_u^T\uPsi\uH_u$. Then we represent $\tilde{\uDelta}$ as follows
\begin{align*}
    \left[ {\begin{array}{cc} \tilde{\Delta}_{11} & (\tilde{\Delta}_{12}^T)^{1 \times (q_n-1)} \\
   \tilde{\Delta}_{12}^{(q_n-1)\times1} & \tilde{\uDelta}_{22}^{(q_n-1)\times(q_n-1)} \\
  \end{array} } \right].
\end{align*}

Now consider the transformation $\tilde{\Delta}_{11}=a,\;\tilde{\Delta}_{12}=ab^T,\;\tilde{\Delta}_{22}=\uC+abb^T$, where $a$ is scalar quantity, $b$ and $\uC$ are $((q_n-1)\times1)$ vector and $((q_n-1)\times(q_n-1))$ matrix respectively. Now we will define our transformation function $h$ as follows,
 \begin{align*}
      h(A_{11},\;A_{12}^{1 \times(q_n-1)},\;\uA_{22}^{(q_n-1)\times(q_n-1)}):=\left[ {\begin{array}{cc}
   A_{11} & A_{11}A_{12} \\
   A_{11}A_{12}^T & \uA_{22}+A_{12}A_{11}^{-1}A_{12}^T \\
  \end{array} } \right].
\end{align*}
Then it follows $h(a,b^T,\uC)=\tilde{\uDelta}$ and Jacobian of transformation for this transformation will be $a^{q_n-1}$ using \cite[Proposition 5.11]{eaton}. Now it is obvious from the formulas of Schur's complement $|\tilde{\uDelta}|=a|\uC|$. Also, if we write $\uS_{Yu}=\left[ {\begin{array}{cc}
   s_{11}^u & s_{12}^u \\
   (s_{12}^{u})^T & \uS_{22} \\
  \end{array} } \right]$ in a similar fashion as $\tilde{\uDelta}$  it is easy to write $\tilde{\uDelta}+\uS_{Yu}=(a+s_{11}^u)\left|\uC^{*}+abb^T-\frac{(s_{12}^u+ab^T)(s_{12}^u+ab^T)^T}{a+s_{11}^u}\right|$, where $\uC^*=\uC+\uS_{22}$. After simple manipulation we can write, $\tilde{\uDelta}+\uS_{Yu}=(a+s_{11}^u)\left|\uC^{**}-\frac{(s_{11}^u)^2}{a+s_{11}^u}\left(b-\frac{s_{12}}{s_{11}^u}\right)\left(b-\frac{s_{12}}{s_{11}^u}\right)^T\right|$, where $\uC^{**}=\uC^*+s_{11}^ubb^T-b(s_{12}^{u})^T-s_{12}^{u}b^T$. Note that $\left(\uC^{**}-\frac{(s_{11}^u)^2}{a+s_{11}^u}\left(b-\frac{s_{12}}{s_{11}^u}\right)\left(b-\frac{s_{12}}{s_{11}^u}\right)^T\right)$ is a positive definite matrix whose determinant is also an increasing function of $a$. 
  
  Using all these facts we can write (\ref{th2:3}) as follows
  \begingroup
\allowdisplaybreaks
  \begin{align}\label{th2:4}
     &\Pi_n(u^T\bar{\uDelta}u > \omega_n| \uYn)=\Pi_n(\tilde{\Delta}_{11}> \omega_n| \uYn)=\Pi_n(a> \omega_n| \uYn)\notag\\
     \leq&\;\frac{\int_{\omega_n}^{\infty}\int_{\mathbb{R}^{(q_n-1)}}\frac{a^{(\nu+\nu_q^*-2)/2+(q_n-1)}}{(a+s_{11}^u)^{(n+\nu+q_n-1)/2}}\exp\left(-a\tr\left(\uPsi^{-1}_u \begin{bmatrix}
    1\\
    b
 \end{bmatrix} \begin{bmatrix}
    1 & b^T  
 \end{bmatrix}\right)/2\right)da\;db}{\int_{0}^{\omega_n}\int_{\mathbb{R}^{(q_n-1)}}\frac{a^{(\nu+\nu_q^*-2)/2+(q_n-1)}}{(a+s_{11}^u)^{(n+\nu+q_n-1)/2}}\exp\left(-a\tr\left(\uPsi^{-1}_u \begin{bmatrix}
    1\\
    b
 \end{bmatrix} \begin{bmatrix}
    1 & b^T  
 \end{bmatrix}\right)/2\right)da\;db}\notag\\
 \times&\;\frac{\int_{\mathbb{P}^{+}_{q_n-1}}|\uC|^{(\nu+\nu_q^*-2)/2}\left|\uC^{**}-\frac{(s_{11}^u)^2}{\omega_n+s_{11}^u}\left(b-\frac{s_{12}}{s_{11}^u}\right)\left(b-\frac{s_{12}}{s_{11}^u}\right)^T\right|^{\frac{(n+\nu+q_n-1)}{2}}\exp\left(\footnotesize-a\tr\left(\uPsi^{-1}_u \begin{bmatrix}
    0 & 0\\
    0 & \uC
 \end{bmatrix}/2\right)\right)dC}{\int_{\mathbb{P}^{+}_{q_n-1}}|\uC|^{(\nu+\nu_q^*-2)/2}\left|\uC^{**}-\frac{(s_{11}^u)^2}{\omega_n+s_{11}^u}\left(b-\frac{s_{12}}{s_{11}^u}\right)\left(b-\frac{s_{12}}{s_{11}^u}\right)^T\right|^{\frac{(n+\nu+q_n-1)}{2}}\exp\left(\footnotesize-a\tr\left(\uPsi^{-1}_u \begin{bmatrix}
    0 & 0\\
    0 & \uC
 \end{bmatrix}/2\right)\right)dC}\notag\\
 \leq&\;\frac{\int_{\omega_n}^{\infty}\frac{a^{(\nu+\nu_q^*-2)/2+(q_n-1)}}{(a+s_{11}^u)^{(n+\nu+q_n-1)/2}}\;da}{\int_9^{\omega_n}\frac{a^{(\nu+\nu_q^*-2)/2+(q_n-1)}}{(a+s_{11}^u)^{(n+\nu+q_n-1)/2}}\;da}\times\frac{\int_{\mathbb{R}^{(q_n-1)}}\exp\left(-\omega_n\tr\left(\uPsi^{-1}_u \begin{bmatrix}
    1\\
    b
 \end{bmatrix} \begin{bmatrix}
    1 & b^T  
 \end{bmatrix}\right)/2\right)\;db}{\int_{\mathbb{R}^{(q_n-1)}}\exp\left(-\omega_n\tr\left(\uPsi^{-1}_u \begin{bmatrix}
    1\\
    b
 \end{bmatrix} \begin{bmatrix}
    1 & b^T  
 \end{bmatrix}\right)/2\right)\;db}\notag\\
 =&\;\frac{\int_{\omega_n}^{\infty}\frac{a^{\frac{(\nu+\nu_q^*-2)}{2}+(q_n-1)}}{(a+s_{11}^u)^{\frac{(n+\nu+q_n-1)}{2}}}\;da}{\int_9^{\omega_n}\frac{a^{\frac{(\nu+\nu_q^*-2)}{2}+(q_n-1)}}{(a+s_{11}^u)^{\frac{(n+\nu+q_n-1)}{2}}}\;da}.
 \end{align}
 \endgroup
 Thus using (\ref{th2:4}), (\ref{th2:2}) will reduced to,
\begin{align}\label{th2:5}
     &\Pi_n(\norm{\bar{\uDelta}} > 2\omega_n | \uYn) \leq\;  9^{q_n} \underset{\norm{u}= 1}{\sup}\Pi_n(u^T\bar{\uDelta}u > \omega_n|\uYn)\notag\\
\leq\;  &9^{q_n}\underset{\norm{u}= 1}{\sup} \frac{\int_{\omega_n}^{\infty}\frac{a^{\frac{(\nu+\nu_q^*-2)}{2}+(q_n-1)}}{(a+s_{11}^u)^{\frac{(n+\nu+q_n-1)}{2}}}\;da}{\int_9^{\omega_n}\frac{a^{\frac{(\nu+\nu_q^*-2)}{2}+(q_n-1)}}{(a+s_{11}^u)^{\frac{(n+\nu+q_n-1)}{2}}}\;da}\leq\; 9\underset{\norm{u}= 1}{\sup} \frac{\int_{\omega_n}^{\infty}\frac{a^{\frac{(\nu+\nu_q^*-2)}{2}+(q_n-1)}}{(a+s_{11}^u)^{\frac{(n+\nu+q_n-1)}{2}}}\;da}{\int_9^{\omega_n}\frac{a^{\frac{(\nu+\nu_q^*-2)}{2}}}{(a+s_{11}^u)^{\frac{(n+\nu+q_n-1)}{2}}}\;da}.
\end{align}
Now recall $s_{11}^u=u^T\uS_{Y}u\in (\lambda_{min}(\uS_Y),\lambda_{max}(\uS_Y))$ as $\norm{u}=1$. Using this fact and Assumption \ref{as7} i.e. $\nu_q^{*}=O(q_n)$ we can similarly proceed as in (\ref{th1:10}-\ref{th1:14}) of Thoerem \ref{th_siw} using the substitution $u_1=\frac{a}{a+s_{11}^u}$ and $k=9$. Hence one can complete this proof using the similar technique mentioned in Theorem \ref{th_siw}. We are not repeating the steps because the calculations are virtually similar. Keeping this note in mind we can say $E_{\mathbb{P}_0}\left[\Pi_n(\norm{\bar{\uDelta}} > 2\omega_n | \uYn)\right]\to 0$ as $n \to \infty$ and proof is done.

\subsection{Proof of Theorem \ref{th_B1} }\label{proof_B1}

\noindent
To begin the proof we first note that,
   \begin{align}\label{th3:1}
       \Pi_n(\norm{\uB-\uB_0} > M\delta_n^{*} | \uYn) \leq &\; \Pi_n\left(\norm{\uB-\tilde{\uB}_n} > M\delta_n^{*}/2 | \uYn\right) +\notag\\ & \Pi_n\left(\norm{\tilde{\uB}_n-\uB_0} > M\delta_n^{*}/2 | \uYn\right),
   \end{align}
where $\tilde{\uB}_n=E[\uB|\uSigma,\uYn]=\uX_{\lambda}^{-1}\uX_n^{T}\uYn$ by Lemma \ref{Lemma_post1} or Lemma \ref{Lemma_post2}. Note that the second term of the right side of the inequality depends only on data. Keeping that in mind we take expectations on both sides over $\uYn$ the above reduces to,
 \begin{align}\label{th3:2}
       E_{\mathbb{P}_0}[\Pi_n(\norm{\uB-\uB_0} > M\delta_n^{*} | \uYn)] \leq &\; E_{\mathbb{P}_0}\left[\Pi_n\left(\norm{\uB-\tilde{\uB}_n} > M\delta_n^{*}/2 | \uYn\right)\right] +\notag\\ &\;\mathbb{P}_0\left(\norm{\tilde{\uB}_n-\uB_0} > M\delta_n^{*}/2\right).
   \end{align}
To show the first term of the right side converges to zero as $n\to\infty$ note that if we define, $\uU_1=\uX_{\lambda}^{1/2}(\uB-\tilde{\uB}_n)\uSigma^{-1/2}$ then $\uU_1|\uSigma,\uYn\sim\mathcal{MN}_{p_n \times q_n}(\uzero,\;\uI_p,\;\uI_q)$.Thus $\norm{\uB-\tilde{\uB}_n}\leq\sqrt{\norm{\uX_{\lambda}^{-1}}}\norm{\uU_1}\sqrt{\norm{\uSigma}}$. Observe the following $\norm{\uX_{\lambda}^{-1}}=\norm{(\uX_n^{T}\uX_n+\lambda \uI_q)^{-1}}\leq\norm{(\uX_n^{T}\uX_n)^{-1}}=\frac{1}{\lambda_{\min}(\uX_n^{T}\uX_n)}\leq\frac{1}{nk_x},$ using Assumption \ref{as5}. Using this fact and using the tower property of conditional expectations we get
\begin{align}\label{th3:3}
    &E_{\mathbb{P}_0}\left[\Pi_n\left(\norm{\uB-\tilde{\uB}_n} > M\delta_n^{*}/2 | \uYn\right)\right]\notag\\
    = &E_{\mathbb{P}_0}\left[E_{\uSigma}\left[\Pi_n\left(\norm{\uB-\tilde{\uB}_n} > M\delta_n^{*}/2 |\uSigma, \uYn\right)|\uYn\right]\right]\notag\\\leq &\;E_{\mathbb{P}_0}\left[E_{\uSigma}\left[\Pi_n\left(\norm{\uU_1} > \frac{M^*\sqrt{n}\delta_n^{*}}{2\sqrt{\beta}} |\uSigma, \uYn\right)|\uYn\right]\right]+E_{\mathbb{P}_0}\left[\Pi_n\left(\norm{\uSigma}>\beta\mid\uYn\right)\right]\notag\\
    \leq &\;E_{\mathbb{P}_0}\left[E_{\uSigma}\left[\Pi_n\left(\norm{\uU_1} > \frac{M^*\sqrt{n}\delta_n^{*}}{2\sqrt{\beta}} |\uSigma, \uYn\right)|\uYn\right]\right]+ E_{\mathbb{P}_0}\left[\Pi_n\left(\norm{\uSigma-\uSigma_0}>\beta/2|\uYn\right)\right]+\notag\\&\uone_{\{\norm{\uSigma_0}>\beta/2\}},
\end{align}
where $\beta$ is a constant chosen such that $\beta>\frac{2}{k_{\sigma}}$ and $M^*=Mk_x$. Then by Assumption \ref{as2} third indicator of the last inequality will be $0$ and using Theorem \ref{th_siw} or Theorem \ref{th_matF}, as appropriate,  $$E_{\mathbb{P}_0}\left[\Pi_n(\norm{\uSigma-\uSigma_0}>\beta/2|\uYn) \right]\to 0$$ as $n \to \infty$. For the first term note that
$\uU_1$ is a matrix whose elements are coming from a i.i.d $N(0,1)$ distribution. So we can think it of as a special case of a matrix with isotropic and independent sub-Gaussian rows with the bounded sub-Gaussian norm. Hence using a similar step as of proving Lemma \ref{subg4}, it follows for sufficiently large $n$,
\begin{align*}
    &E_{\mathbb{P}_0}\left[E_{\uSigma}\left[\Pi_n\left(\norm{\uU_1} > \frac{M^*\sqrt{n}\delta_n^{*}}{2\sqrt{\beta}} |\uSigma, \uYn\right)|\uYn\right]\right]\\\leq&\;E_{\mathbb{P}_0}\left[E_{\uSigma}\left[\Pi_n\left(\norm{\uU_1} > \frac{M^*\max{(\sqrt{p_n},\sqrt{q_n})}}{2\sqrt{\beta}} |\uSigma, \uYn\right)|\uYn\right]\right]\\\leq\;&\exp{\left(-\frac{c_{\sigma_{0}\max{(\sqrt{p_n},\sqrt{q_n})}}}2\right)}\to 0,\;
\end{align*}
as $n\to \infty$. It now remains to be shown that the second term of (\ref{th3:2}) converges to zero. Note that
\begin{align}\label{th3:4}
    \mathbb{P}_0\left(\norm{\tilde{\uB}_n-\uB_0} > M\delta_n^{*}/2\right)\leq \mathbb{P}_0\left(\norm{\tilde{\uB}_n-\BLS} > M\delta_n^{*}/4\right)+\mathbb{P}_0\left(\norm{\BLS-\uB_0} > M\delta_n^{*}/4\right)
\end{align}
For the second term of the right side of the inequality it follows from Lemma \ref{subg4},
\begin{align}\label{th3:5}
    \mathbb{P}_0\left(\norm{\BLS-\uB_0} > M\delta_n^{*}/4\right) \to 0 \;\textit{as}\; n \to \infty. 
\end{align}
For the first term of (\ref{th3:4}) note that $\tilde{\uB}_n-\BLS=(\uX_{\lambda}^{-1}-(\uX_n^{T}\uX_n)^{-1})\uX_n^{T}\uYn$. Now by Woodbury's identity it follows
$\uX_{\lambda}^{-1}-(\uX_n^{T}\uX_n)^{-1}=-(\uX_n^{T}\uX_n)^{-1}\uX_{\lambda}^{*-1}(\uX_n^{T}\uX_n)^{-1},$ where $\uX_{\lambda}^*=((\uX_n^{T}\uX_n)^{-1}+\lambda^{-1}\uI_p)$. Then,
 $\BLS-\tilde{\uB}_n=(\uX_n^{T}\uX_n)^{-1}\uX_{\lambda}^{*-1}\BLS.$
 Now since $\norm{(\uX_n^{T}\uX_n)^{-1}}\norm{\uX_{\lambda}^{*-1}}\leq\frac{\lambda_{max}(\uX_n^{T}\uX_n)}{\lambda_{min}(\uX_n^{T}\uX_n)}\leq\frac{1}{k_x^2}$ by Assumption \ref{as5}. Then it follows that
\begin{align}\label{th3:6}
    &\mathbb{P}_0\left(\norm{\tilde{\uB}_n-\BLS} > M\delta_n^{*}/4\right)\leq \mathbb{P}_0\left(\norm{(\uX_n^{T}\uX_n)^{-1}}\norm{\uX_{\lambda}^{*-1}}\norm{\BLS}>M\delta_n^{*}/4\right)
    \notag\\\leq &\; \mathbb{P}_0\left(\norm{(\uX_n^{T}\uX_n)^{-1}}\norm{\uX_{\lambda}^{*-1}}\norm{\BLS-\uB_0}>M\delta_n^{*}/8\right)+\notag\\&\;\uone_{\left\{ \norm{{\left(\uX_n^{T}\uX_n\right)^{-1}}}\norm{{\uX_{\lambda}^{*-1}}}{\norm{\uB_{0}}>\frac{M\delta_n^{*}}{8}}\right\}}\notag\\ \leq & \;\mathbb{P}_0 \left(\norm{\BLS-\uB_0}>\frac{M\delta_n^{*}k_x^2}{8}\right)+ \uone_{\left\{\norm{(\uX_n^{T}\uX_n)^{-1}}\norm{{\uX_{\lambda}^{*-1}}}\norm{{\uB_{0}}}> \frac{M\delta_n^{*}}{8}\right\}}.
\end{align}
The first term of the right side of (\ref{th3:6}) $\to 0$ as $n \to \infty$ similarly as in (\ref{th3:5}). For the second term note that $\norm{(\uX_n^{T}\uX_n)^{-1}}\norm{\uX_{\lambda}^{*-1}}\leq\frac{(\lambda_{min}((\uX_n^{T}\uX_n)^{-1})+\lambda^{-1})^{-1}}{\lambda_{min}(\uX_n^{T}\uX_n)}\leq\frac{(k_x+n\lambda^{-1})^{-1}}{k_x}$. It follows that for sufficiently large n, 
\begin{align*}
    \uone_{\left\{\norm{(\uX_n^{T}\uX_n)^{-1}}\norm{\uX_{\lambda}^{*-1}}\norm{\uB_{0}}>M\delta_n^{*}/8\right\}}\leq\uone_{\left\{\lambda^{-1}<\frac{8\norm{\uB_{0}}}{Mn\delta_n^{*}k_x}-\frac{k_x}{n}\right\}}\leq \uone_{\left\{\lambda^{-1}<\frac{\lambda_0 \norm{\uB_{0}}}{\max{(\sqrt{np_n},\sqrt{nq_n})}}-\frac{k_x}{n}\right\}},
\end{align*}
 where $\lambda_0=\frac{8}{Mk_x}$. One can make the last indicator term $0$ for sufficiently large $n$ using Assumption \href{a}{6}. This completes the proof.

\subsection{Proof of Theorem \ref{th_B2} }\label{proof_B2}

\noindent
To initiate the proof, it is important to observe that:
   \begin{align}\label{th4:1}
       \Pi_n(\fnorm{\uB-\uB_0} > M\delta_n^{**} | \uYn) \leq &\; \Pi_n\left(\fnorm{\uB-\tilde{\uB}_n} > M\delta_n^{**}/2 | \uYn\right) +\notag\\ & \Pi_n\left(\fnorm{\tilde{\uB}_n-\uB_0} > M\delta_n^{**}/2 | \uYn\right),
   \end{align}
where $\tilde{\uB}n=E[\uB|\uSigma,\uYn]=\uX_{\lambda}^{-1}\uX_n^{T}\uYn$ according to either Lemma \ref{Lemma_post1} or Lemma \ref{Lemma_post2}. Notably, the second term on the right side of the inequality is solely dependent on the data. Considering this, we proceed by taking expectations on both sides over $\uYn$, which simplifies to:
 \begin{align}\label{th4:2}
       E_{\mathbb{P}_0}[\Pi_n(\fnorm{\uB-\uB_0} > M\delta_n^{**} | \uYn)] \leq &\; E_{\mathbb{P}_0}\left[\Pi_n\left(\fnorm{\uB-\tilde{\uB}_n} > M\delta_n^{**}/2 | \uYn\right)\right] +\notag\\ &\;\mathbb{P}_0\left(\fnorm{\tilde{\uB}_n-\uB_0} > M\delta_n^{**}/2\right).
   \end{align}
To demonstrate the convergence of the first term on the right side to zero as $n\to\infty$, consider defining $\uU_1=\uX_{\lambda}^{1/2}(\uB-\tilde{\uB}n)\uSigma^{-1/2}$. It follows that $\uU_1|\uSigma,\uYn\sim\mathcal{MN}{p_n \times q_n}(\uzero,;\uI_p,;\uI_q)$. Consequently, $\fnorm{\uB-\tilde{\uB}n}\leq\sqrt{\norm{\uX{\lambda}^{-1}}}\fnorm{\uU_1}\sqrt{\norm{\uSigma}}$. Considering $\norm{\uX_{\lambda}^{-1}}$, we have $\norm{\uX_{\lambda}^{-1}}\leq\frac{1}{nk_x}$ using Assumption \ref{as5}. By utilizing the tower property of conditional expectations, the first term on the right side can be expressed as:
\begin{align}\label{th4:3}
    &E_{\mathbb{P}_0}\left[\Pi_n\left(\fnorm{\uB-\tilde{\uB}_n} > M\delta_n^{**}/2 | \uYn\right)\right]\notag\\
    = &E_{\mathbb{P}_0}\left[E_{\uSigma}\left[\Pi_n\left(\fnorm{\uB-\tilde{\uB}_n} > M\delta_n^{**}/2 |\uSigma, \uYn\right)|\uYn\right]\right]\notag\\\leq &\;E_{\mathbb{P}_0}\left[E_{\uSigma}\left[\Pi_n\left(\fnorm{\uU_1} > \frac{M^*\sqrt{n}\delta_n^{**}}{2\sqrt{\beta}} |\uSigma, \uYn\right)|\uYn\right]\right]+E_{\mathbb{P}_0}\left[\Pi_n\left(\norm{\uSigma}>\beta\mid\uYn\right)\right]\notag\\
    \leq &\;E_{\mathbb{P}_0}\left[E_{\uSigma}\left[\Pi_n\left(\fnorm{\uU_1} > \frac{M^*\sqrt{n}\delta_n^{**}}{2\sqrt{\beta}} |\uSigma, \uYn\right)|\uYn\right]\right]+ E_{\mathbb{P}_0}\left[\Pi_n\left(\norm{\uSigma-\uSigma_0}>\beta/2|\uYn\right)\right]+\notag\\&\uone_{\{\norm{\uSigma_0}>\beta/2\}},
\end{align}
where $\beta$ is a constant chosen such that $\beta>\frac{2}{k_{\sigma}}$ and $M^*=Mk_x$. Using Assumption \ref{as2}, the third indicator in the last inequality becomes $0$, and by applying Theorem \ref{th_siw} or Theorem \ref{th_matF} appropriately, it follows that $E_{\mathbb{P}_0}[\Pi_n(\norm{\uSigma-\uSigma_0}>\beta/2|\uYn)] \to 0$ as $n \to \infty$. For the first term, consider that $\uU_1$ is a matrix with elements from an i.i.d $N(0,1)$ distribution. Viewing it as a special case of a matrix with isotropic and independent sub-Gaussian rows and a bounded sub-Gaussian norm, it follows, by a similar approach as in proving Lemma \ref{subg4}, that:
\begin{align*}
    &E_{\mathbb{P}_0}\left[E_{\uSigma}\left[\Pi_n\left(\fnorm{\uU_1} > \frac{M^*\sqrt{n}\delta_n^{**}}{2\sqrt{\beta}} |\uSigma, \uYn\right)|\uYn\right]\right]\\\leq&\;E_{\mathbb{P}_0}\left[E_{\uSigma}\left[\Pi_n\left(\fnorm{\uU_1} > \frac{M^*\sqrt{p_nq_n}}{2\sqrt{\beta}} |\uSigma, \uYn\right)|\uYn\right]\right]\\\leq\;&\exp{\left(-c_2p_nq_n\right)}\to 0,\;
\end{align*}
as $n\to \infty$. It remains to be shown that the second term of (\ref{th4:2}) converges to zero. Note that
\begin{align}\label{th4:4}
    \mathbb{P}_0\left(\fnorm{\tilde{\uB}_n-\uB_0} > M\delta_n^{**}/2\right)\leq \mathbb{P}_0\left(\fnorm{\tilde{\uB}_n-\BLS} > M\delta_n^{**}/4\right)+\mathbb{P}_0\left(\fnorm{\BLS-\uB_0} > M\delta_n^{**}/4\right)
\end{align}
For the second term of the right side of the inequality it follows from Lemma \ref{subg4},
\begin{align}\label{th4:5}
    \mathbb{P}_0\left(\fnorm{\BLS-\uB_0} > M\delta_n^{**}/4\right) \to 0 \;\textit{as}\; n \to \infty. 
\end{align}
For the first term in (\ref{th4:4}), we can follow a similar approach as in the proof of Theorem \ref{th_B1}. By employing (\ref{th4:5}) and leveraging Assumption \ref{as8}, we can demonstrate that $\mathbb{P}_0\left(\fnorm{\tilde{\uB}_n-\BLS} > M\delta_n^{**}/4\right)\to0$ as $n\to\infty$. This concludes the proof.

\section{Discussion}\label{sec10}
This article describes and examines the theoretical properties of posterior distributions corresponding to DSIW priors and matrix-$F$ priors for the error covariance matrix in a multi-response regression setting. We establish posterior consistency for $\uSigma$ with the standard posterior contraction rate (Theorem \ref{th_siw} and Theorem \ref{th_matF}) under mild regularity assumptions on the number of variables, the true underlying covariance matrix, the design matrix, and relevant hyperparameters. Although our main parameter of interest is $\uSigma$, we also establish a posterior contraction rate for coefficient matrix $\uB$ (Theorem \ref{th_B1}, \ref{th_B2}) under standard regularity assumptions. The results related to the inconsistency of posterior distribution in Section \ref{sec7} demonstrate that extension of our results to ultra-high dimensional settings is not possible without imposing additional low-dimensional structure. In Section \ref{sec8}, we established the geometric ergodicity of Gibbs samplers for both DSIW and matrix-$F$ prior distributions. This provides key convergence guarantees and valid estimates of standard error for various mcmc-based approximations of posterior quantities.
The simulation study strongly validates our technical results (Section \ref{sec9}).

Most DSIW priors in the literature use a fixed value for the degrees of freedom hyperparameter $\nu$. However, some authors (see for example \cite{daniels, feron}) suggest using a prior for the degrees of freedom $\nu$ of DSIW prior (we assume $\nu$ to be fixed in our analysis). Our analysis can be easily extended to include settings where the prior distribution on degrees of freedom $\nu$ has support $(0,\;m_n)$ with $m_n = o(n)$. However, non-trivial additional work is needed when the prior for $\nu$ has unbounded support, and this will be investigated as part of future work. Another future direction is to develop sparse or low-rank versions of the DSIW and matrix-F priors and investigate their empirical and theoretical properties.

%
%

\begin{acks}[Acknowledgments]
The authors express their gratitude to Andrew Gelman and Matt Wand for their supportive and constructive feedback, which significantly enhanced the quality of the paper. We would also like to thank three anonymous referees for their valuable comments that improved the quality of the paper.
\end{acks}
%

\begin{supplement}
\stitle{Supplement to “Posterior consistency in multi-response regression models with non-informative priors for the error covariance matrix in growing dimensions”}
\sdescription{The supplement (\cite{supp}) provides the remaining proofs.}
\end{supplement}


\bibliographystyle{imsart-nameyear} 
\bibliography{bibfile}       


\end{document}